\documentclass[a4paper,12pt,final]{amsart}

\usepackage{amssymb}
\usepackage{latexsym}
\usepackage{amsmath}
\usepackage{mathrsfs}
\usepackage{mathtools}
\usepackage[all]{xy}
\usepackage{enumitem}
\usepackage{graphicx}
\usepackage{pdflscape}

\usepackage{adjustbox}

\usepackage{xcolor}
\usepackage{hyperref}
\hypersetup{
    colorlinks,
    linkcolor={red!50!black},
    citecolor={blue!50!black},
    urlcolor={blue!80!black}
}

\newcommand{\harxiv}[1]{\href{http://arxiv.org/abs/#1}{\texttt{arXiv:#1}}}
\newcommand{\hyref}[2]{ \hyperref[#2]{#1~\ref*{#2}} }

\theoremstyle{boldhead}

\theoremstyle{plain}
\newtheorem{theorem}{Theorem}[section]
\newtheorem{lemma}[theorem]{Lemma}
\newtheorem{corollary}[theorem]{Corollary}
\newtheorem{proposition}[theorem]{Proposition}
\newtheorem{question}[theorem]{Question}
\newtheorem{introtheorem}{Theorem}

\theoremstyle{definition}
\newtheorem{remark}[theorem]{Remark}
\newtheorem{example}[theorem]{Example}
\newtheorem{definition}[theorem]{Definition} 
\newtheorem{setup}[theorem]{Setup}

\newcommand{\sA}{\mathsf{A}}
\newcommand{\sB}{\mathsf{B}}
\newcommand{\sD}{\mathsf{D}}
\newcommand{\sH}{\mathsf{H}}
\newcommand{\sI}{\mathsf{I}}
\newcommand{\sK}{\mathsf{K}}
\newcommand{\sL}{\mathsf{L}}
\newcommand{\sP}{\mathsf{P}}
\newcommand{\sS}{\mathsf{S}}
\newcommand{\sU}{\mathsf{U}}
\newcommand{\sV}{\mathsf{V}}
\newcommand{\sW}{\mathsf{W}}
\newcommand{\sX}{\mathsf{X}}
\newcommand{\sY}{\mathsf{Y}}
\newcommand{\sZ}{\mathsf{Z}}

\newcommand{\Db}{\mathsf{D}^b}
\newcommand{\Kb}{\mathsf{K}^b}
\newcommand{\Dp}{\mathsf{D}^p}
\newcommand{\Dsg}{\mathsf{D}_{sg}}
\newcommand{\Zsg}{\mathsf{Z}_{sg}}

\newcommand{\cF}{\mathcal{F}}
\newcommand{\cP}{\mathcal{P}}
\newcommand{\cT}{\mathcal{T}}

\newcommand{\bH}{\mathbb{H}}
\newcommand{\bN}{\mathbb{N}}
\newcommand{\bR}{\mathbb{R}}
\newcommand{\bC}{\mathbb{C}}
\newcommand{\bP}{\mathbb{P}}
\newcommand{\bS}{\mathbb{S}}
\newcommand{\bZ}{\mathbb{Z}}

\newcommand{\Hminus}{\bH}

\renewcommand{\geq}{\geqslant}
\renewcommand{\leq}{\leqslant}
\renewcommand{\phi}{\varphi}
\renewcommand{\epsilon}{\varepsilon}

\DeclareMathOperator{\Hom}{\mathsf{Hom}}
\DeclareMathOperator{\Ext}{\mathsf{Ext}}

\newcommand{\add}[1]{\mathsf{add}(#1)}
\renewcommand{\mod}[1]{\mathsf{mod}(#1)}
\newcommand{\stmod}[1]{\underline{\mathsf{mod}}(#1)}
\newcommand{\repb}[1]{\mathsf{rep}^b(#1)}
\newcommand{\proj}[1]{\mathsf{proj}(#1)}
\newcommand{\coh}[1]{\mathsf{coh}(#1)}
\newcommand{\inj}[1]{\mathsf{inj}(#1)}
\newcommand{\nil}[1]{\mathsf{nilrep}(#1)}

\DeclareMathOperator{\coker}{\mathrm{coker}}

\newcommand{\simples}[1]{\mathsf{Simp}(#1)}

\newcommand{\smfull}{\sU}
\newcommand{\smsub}{\sS}

\newcommand{\kk}{{\mathbf{k}}}
\newcommand{\dd}{{\mathbf{d}}}

\newcommand{\SSS}{\mathbb{S}}

\newcommand{\extn}[1]{\langle #1 \rangle}

\newcommand{\orth}{^\perp}
\newcommand{\orthw}{^{\perp_w}}

\DeclareMathOperator{\Susp}{\mathsf{susp}}
\DeclareMathOperator{\Cosusp}{\mathsf{cosusp}}
\newcommand{\susp}[2]{\Susp_{#1}(#2)}
\newcommand{\cosusp}[2]{\Cosusp_{#1}(#2)}

\newcommand{\cocone}[1]{\mathsf{cocone}(#1)}

\newcommand{\thick}[2]{\mathsf{thick}_{#1}(#2)}

\renewcommand{\tt}{\mathbf{t}}

\newcommand{\lsms}[2]{\mathsf{L}_{#1}(#2)}
\newcommand{\lsmsn}[3]{\mathsf{L}^{#1}_{#2}(#3)}
\newcommand{\rsms}[2]{\mathsf{R}_{#1}(#2)}
\newcommand{\rsmsn}[3]{\mathsf{R}^{#1}_{#2}(#3)}
\newcommand{\shift}[1]{{\langle #1 \rangle}}  

\newcommand{\lhrs}[2]{\mathsf{Ltilt}_{#1}(#2)}

\newcommand{\rhrs}[2]{\mathsf{Rtilt}_{#1}(#2)}

\newcommand{\too}{\longrightarrow}
\newcommand{\rightlabel}[1]{\stackrel{#1}{\longrightarrow}}
\newcommand{\bij}{\stackrel{1-1}{\longleftrightarrow}}
\newcommand{\into}{\hookrightarrow}
\newcommand{\onto}{\twoheadrightarrow}
   
\entrymodifiers={+!!<0pt,\fontdimen22\textfont2>}

\newcommand{\coloneqq}{\mathrel{\mathop:}=}

\usepackage{tikz}

\usepackage{tikz}
\usetikzlibrary{decorations.pathmorphing}
\usetikzlibrary{decorations.pathreplacing}
\usetikzlibrary{positioning}
\usetikzlibrary{calc}
\usetikzlibrary{backgrounds}
\usetikzlibrary{shapes.misc}
\usetikzlibrary{patterns}
\usetikzlibrary{math}

\usepackage{tikz,tikz-cd}
\tikzcdset{scale cd/.style={every label/.append style={scale=#1},
    cells={nodes={scale=#1}}}}
    
\makeatletter
\@namedef{subjclassname@2020}{%
  \textup{2020} Mathematics Subject Classification}
\makeatother

\title[Simple tilts and mutation]{Simple tilts of length hearts and simple-minded mutation}
\author{Nathan Broomhead, Raquel Coelho Sim\~oes, David Pauksztello, Jon Woolf}

\makeatletter
\@namedef{subjclassname@2020}{%
  \textup{2020} Mathematics Subject Classification}
\makeatother

\begin{document}

\begin{abstract}
We characterise when a simple Happel--Reiten--Smal\o\ tilt of a length heart is again a length heart in terms of approximation theory and the existence of a stability condition with a phase gap.
We apply simple-minded reduction to provide a sufficient condition for infinite iterability of simple-minded mutation/simple tilting. 
We use simple-minded mutation pairs to provide a common framework to show that mutation of simple-minded collections (resp. $w$-simple-minded systems, for $w \geq 1$) gives simple-minded collections (resp. $w$-simple-minded systems) under mild conditions, in the process providing a unified proof of results of Alex Dugas \cite{Dugas} and Peter J\o rgensen \cite{J22}.
Finally, we show that under mild conditions, mutation of simple-minded collections is compatible with mutation of $w$-simple-minded systems via a singularity category construction due to Haibo Jin \cite{Jin}.
\end{abstract}



\maketitle

{\small
  \tableofcontents
}

\addtocontents{toc}{\protect\setcounter{tocdepth}{1}}

\section*{Introduction}

\noindent

Homological algebra provides a common framework for many branches of mathematics and as such, the language of triangulated categories and abelian categories is widely used to study representation theory, algebraic geometry, symplectic geometry and algebraic topology. Classic tilting theory describes derived equivalences in terms of tilting objects. 
Indeed, Alexander Beilinson's famous derived equivalence $\Db(\coh{\bP^1}) \simeq \Db(\kk \widetilde{A}_1)$ in \cite{Beilinson78}, where $\kk \widetilde{A}_1$ is the Kronecker algebra, can be considered the first theorem of tilting theory, providing an unexpected and deep connection between representation theory and algebraic geometry.
In this article we consider two generalisations of tilting theory and their interaction.

The first generalisation is \emph{Happel--Reiten--Smal\o\ (HRS) tilting}. 
It is formulated in the language of t-structures and torsion pairs. Let $\sD$ be a triangulated category with shift functor $[1]\colon \sD \to \sD$. 
A \emph{t-structure} $(\sX,\sY)$ in $\sD$ is a pair of subcategories giving rise to a cohomology theory on $\sD$ taking values in an abelian category $\sH = \sX \cap \sY[1]$ called the \emph{heart}.
A \emph{torsion pair} $(\cT,\cF)$ in $\sH$ is a framework for abelian categories that abstracts the properties of torsion and torsionfree abelian groups in the category of finitely generated abelian groups.
For a t-structure $(\sX,\sY)$ and a torsion pair $\tt = (\cT,\cF)$ in its heart $\sH$, the \emph{right HRS tilt of $\sH$ at $\tt$} is the new t-structure
\[
(\sX * \cT[-1], (\cF * \sY)[-1]) 
\quad \text{with new heart} \quad
\sK = \cF * \cT[-1].
\]
When $\cT$ is the extension closure of a set of simple objects of $\sH$, we call $\sK$ a \emph{right simple tilt} of $\sH$.
HRS tilting has two main applications --- it provides:
\begin{itemize}
\item  a `mutation theory' for t-structures, see \cite{HRS}; and,
\item a method for doing tilting theory without tilting objects, i.e. constructing derived equivalences when no tilting objects are available, see \cite{CHZ19}.
\end{itemize}
The first application is our principal motivation and the aspect we study in this article. It is in this context that it has been used extensively in the study of Bridgeland stability conditions, see \cite{Bridgeland08,CHQ23,PSZ18,QW18,Woolf10}.

In the theory of stability conditions,  \emph{length categories}, i.e.\ abelian categories in which each object is both noetherian and artinian, are of central importance for three main reasons. 
First, categories of semistable objects are naturally length categories. 
Second, any stability function on a length category automatically satisfies the Harder--Narasimhan property.
Third, the geometry of the stability manifold corresponding to crossing a type II wall associated to a simple HRS tilt of a length category which is again length is particularly well behaved \cite{Woolf10}. 
The second observation, in particular, makes it much easier to define stability conditions on length categories.

It is natural, therefore, to ask when a simple HRS tilt of a length category is again length. 
Unfortunately, there is an explicit counterexample in \cite{Koenig-Yang} showing this does not always happen, see Example~\ref{ex:koenig-yang}. Our first main result characterises when it does happen.

\begin{introtheorem}[Theorem~\ref{thm:length-simple-tilts}] \label{intro:length-simple-tilts}
Let $\sH$ be a bounded length heart in $\sD$ whose set of isoclasses of simple objects $\simples{\sH}$ is finite. 
Suppose $(\cT,\cF)$ is a torsion pair such that $\cT = \extn{\smsub}$ is the extension closure of $\smsub \subseteq \simples{\sH}$. Let $\sK = \cF * \cT[-1]$ be the right simple tilt at $\cT$. 
Then the following conditions are equivalent.
\begin{enumerate}[label=(\roman*)]
\item $\sK$ is a length heart. 
\item Each object in $(\simples{\sH} \setminus \smsub)[1]$ admits a right $\extn{\smsub}$-approximation. \label{approx}
\item There exists a stability condition $\sigma = (Z,\cP)$ with $\cP(0,1] = \sH$, $\cP(1)=\extn{\smsub}$ and $\cP(0,\phi] = \smsub\orth\cap \sH$ for some $0<\phi<1$.
\end{enumerate}
\end{introtheorem}

There is an obvious dual statement for left simple tilts.
The second statement above involves approximation theory, which is a key ingredient in \emph{mutation} \cite{AI}; for precise definitions see Section~\ref{sec:background}.
Closing the set of tilting objects under the mutation operation gives rise to \emph{silting objects}. Together with \emph{cluster-tilting objects}, they can be considered `projective-minded' objects, see \cite{CKL15}. 
Cluster-tilting objects can be thought of as the `shadow' of silting objects in Calabi--Yau triangulated categories via a singularity category type construction, see \cite{IY}.
For projective-minded objects there are extensive theories of mutation which lead to rich combinatorial structure in the associated representation theory and homological algebra.

However, for length categories, the natural kinds of objects are `simple-minded', i.e. objects that have the homological properties of simple objects.
This leads us to our second generalisation of tilting theory.
\emph{Simple-minded collections (SMCs)}, defined in \cite{Koenig-Yang}, are collections of simple objects in length hearts.
Their shadow in Calabi--Yau triangulated categories are \emph{$w$-simple-minded systems ($w$-SMSs)}, defined in \cite{Koenig-Liu} for $w = 1$ and in \cite{CS15} for $w >1$. 
SMSs axiomatise the image of the nonprojective simple modules in the stable module category of a selfinjective algebra.

In \cite{Koenig-Yang}, Steffen Koenig and Dong Yang show that mutation of SMCs is always defined and again an SMC when $\sD = \Db(A)$ for a finite-dimensional algebra $A$. As a consequence, for $\Db(A)$, mutation of SMCs can be iterated indefinitely.
For $w$-SMSs, and in the context of $(-w)$-Calabi-Yau triangulated categories, mutation is shown to be defined and again a $w$-SMS in \cite{Dugas} for $w = 1$ and in \cite{J22} for $w > 1$. 
There has been recent work in \cite{Sun-Zhang} investigating compatibility of mutation of SMCs with recollements of triangulated categories and in \cite{Chen} investigating mutations of SMCs in derived categories of tube categories.

Our second main result is an application of simple-minded mutation pairs in \cite{CSP20} and simple-minded reduction in \cite{CSP20,CSPP,Jin}. It provides necessary and sufficient conditions for mutations of SMCs to be infinitely iterable and an alternative proof that mutations of $w$-SMSs are again $w$-SMS for $w \geq 1$ in a more general context, unifying the proof for $w$-SMSs and SMCs. 

\begin{introtheorem}[Theorem~\ref{thm:mut-SM-is-SM}] \label{intro:mut-SM-is-SM}
Let $\sD$ be a Hom-finite, Krull--Schmidt triangulated category.
Suppose $\smfull$ is an SMC or $w$-SMS and $\smsub \subseteq \smfull$. 
Then the mutation of $\smfull$ at $\smsub$ is infinitely iterable if and only if 
\begin{enumerate}
\item for $\smfull$ an SMC, $\extn{\smsub}$ is covariantly finite in ${}\orth (\smsub[{<}0])$ and contravariantly finite in $(\smsub[{>}0])\orth$;
\item for $\smfull$ a $w$-SMS, $\extn{\smsub}$ is functorially finite and $\smsub$ is invariant under the shift of the Serre functor, $\SSS[w]$.
\end{enumerate}
\end{introtheorem}

These conditions are precisely what is required to be able to perform simple-minded reduction. This provides a simple conceptual explanation of why mutation works, the process is `reduce, shift, lift':
\[
\begin{tikzcd}
\left\{ \parbox{3.75cm}{\centering SMCs or $w$-SMSs in $\sD$ containing $\smsub$} \right\} \arrow[d,"\text{reduction}"'] \arrow[rr, leftrightarrow, "\text{right mutation}", "\text{left mutation}"'] & & 
\left\{ \parbox{3.75cm}{\centering SMCs or $w$-SMSs in $\sD$ containing $\smsub$} \right\} \arrow[d,"\text{reduction}"] \\
\{\text{SMCs or $w$-SMSs in $\sZ$}\} \arrow[rr, leftrightarrow, "\shift{1}", "\shift{-1}"']  & & 
\{\text{SMCs or $w$-SMSs in $\sZ$}\}
\end{tikzcd}
\]
Here $\sZ$ is the simple-minded reduction of $\sD$ at $\smsub$ and $\shift{1}$ denotes the shift functor in $\sZ$, see \cite{CSP20,CSPP,Jin} for details.

An SMC corresponds to a bounded t-structure $(\sX,\sY)$ with a length heart. By the Koenig--Yang correspondences in \cite{Koenig-Yang}, when $\sD = \Db(A)$ for a finite-dimensional algebra $A$, any such t-structure admits a left adjacent co-t-structure $({}\orth \sX, \sX)$ and a right adjacent co-t-structure $(\sY, \sY\orth)$. 
We call an SMC whose associated t-structure admits both adjacent co-t-structures \emph{bisilting}. Our third main theorem says that the property of being a finite bisilting SMC is preserved under mutation and provides a conceptual explanation of why mutation is compatible with the Koenig--Yang correspondences in \cite{Koenig-Yang}.

\begin{introtheorem}[Theorem~\ref{thm:bisilting-smc-mutation}] \label{intro:bisilting-smc-mutation}
Let $\sD$ be a Hom-finite, Krull--Schmidt triangulated category. 
The right mutation of a finite bisilting SMC at any subset is also a finite bisilting SMC.
\end{introtheorem}

We note that a similar conceptual consideration of the Koenig--Yang correspondences without reference to mutation has recently been considered in \cite{Bonfert}.

Our final main result is the application of Theorems~\ref{intro:mut-SM-is-SM} and \ref{intro:bisilting-smc-mutation} to obtain compatibility between SMC and $w$-SMS mutation in the setting of a $(1-w)$-Calabi--Yau triple $(\sD,\Dp,\smfull)$ consisting of a triangulated category $\sD$ with thick subcategory $\Dp$ and finite bisilting SMC $\smfull$, see \cite{Jin}. The resulting Verdier quotient $\Dsg = \sD/\Dp$ is a $(-w)$-Calabi--Yau category in which $\smfull$ is a $w$-SMS \cite{Jin}. The prototypical example is $(\Db(A), \Kb(\proj{A}), \simples{A})$ for a finite-dimensional symmetric algebra $A$. Here the Verdier quotient is the singularity category $\Dsg \simeq \stmod{A}$.

\begin{introtheorem}[Theorem~\ref{thm:compatibility}] \label{intro:compatibility}
Let $(\sD, \Dp, \smfull)$ be a $(1-w)$-CY-triple in which $\smfull$ is finite. Suppose $\smsub$ is an $\infty$-orthogonal collection (see Definition~\ref{def:orth}) satisfying mild conditions. Then the following diagram is commutative:
\[
\begin{tikzcd}
\{\text{bisilting SMCs in $\sD$ containing $\smsub$}\} \arrow[rr, "\text{right mutation}"] \ar[d, "\text{Verdier localisation}"'] & &
\{\text{bisilting SMCs in $\sD$ containing $\smsub$}\} \arrow[d, "\text{Verdier localisation}"] \\
\{\text{$w$-SMSs in $\Dsg$ containing $\smsub$}\} \arrow[rr, "\text{right mutation}"']  & &
\{\text{$w$-SMSs in $\Dsg$ containing $\smsub$}\}.
\end{tikzcd}
\]
\end{introtheorem}

Finally we note that preliminary versions of these results were announced in \cite{Pauk23}.

\subsection*{Acknowledgments}
It is a pleasure to thank David Ploog for useful conversations and his comments.
We extend our gratitude to the referee for a careful reading of the paper and helpful comments.
This project has been supported by the European Union's Horizon 2020 research and innovation programme through the Marie Sk\l odowska-Curie Individual Fellowship grant 838706 of the second author and by EPSRC grant no. EP/V050524/1 of the third author.

\section{Background} \label{sec:background}

\noindent
Throughout, $\kk$ denotes a field and $\sD$ will be a Hom-finite, Krull--Schmidt, $\kk$-linear triangulated category with shift functor $[1]\colon \sD \to \sD$. 
For subcategories $\sX$ and $\sY$ of $\sD$ we write
\[
\sX * \sY \coloneqq \{ d \in \sD \mid \text{there is a triangle } x \to d \to y \to x[1] \text{ with } x \in \sX \text{ and } y \in \sY\}.
\]
We say that a subcategory $\sX$ is \emph{extension-closed} if $\sX * \sX = \sX$. 
For a subcategory or collection of objects $\sX$ of $\sD$, we write $\extn{\sX}$ for the \emph{extension closure of $\sX$ in $\sD$}. The \emph{right} and \emph{left perpendicular categories} of $\sX$ are defined as follows:
\[
\sX^\perp = \{d \in \sD \mid \Hom(\sX,d) = 0 \} \text{ and } \,^\perp \sX = \{d \in \sD \mid \Hom(d,\sX)=0 \},
\]
where we employ the shorthand $\Hom(\sX,d) = 0$ to mean $\Hom(x, d) = 0$ for each object $x$ in $\sX$.

A \emph{Serre functor} on $\sD$ is an autoequivalence $\bS\colon \sD \to \sD$ such that there is an isomorphism,
$\Hom(x, y) \simeq D \Hom(y, \bS x)$, which is natural in $x$ and $y$, where $D$ denotes the standard vector space duality. 
For $w \in \bZ$, a triangulated category $\sD$ is called \emph{$w$-Calabi-Yau} if it has a Serre functor $\bS$ and there is a natural isomorphism $\bS \simeq [w]$.

\subsection{Approximations}

In order to define mutations later we need some approximation theory. Let $\sX$ be a full additive subcategory of $\sD$ and $d \in \sD$.
\begin{itemize}
\item A morphism $\alpha \colon x \to d$ with $x \in \sX$ is called a \emph{right $\sX$-approximation} if the induced map $\Hom_\sD(\sX, \alpha) \colon \Hom_{\sD}(\sX, x) \onto \Hom_{\sD}(\sX, d)$ is a surjection. 
\item The morphism $\alpha$ is \emph{right minimal} if any $\beta \colon x \to x$ satisfying $\alpha \beta = \alpha$ is an automorphism. 
\item A morphism $\alpha \colon x \to d$ is a \emph{minimal right $\sX$-approximation} if it is both a right $\sX$-approximation and right minimal.
\item The subcategory $\sX$ is \emph{contravariantly finite} in $\sD$ if every object in $\sD$ admits a right $\sX$-approximation. Since $\sD$ is Hom-finite and Krull--Schmidt, every right $\sX$-approximation contains a minimal right $\sX$-approximation as a summand, see e.g.~\cite{Auslander-Smalo}.
\end{itemize}
There are dual notions of \emph{(minimal) left $\sX$-approximation} and \emph{covariantly finite} subcategories.
The subcategory $\sX$ is \emph{functorially finite} if it is both contravariantly finite and covariantly finite.

\subsection{Orthogonal collections, torsion pairs and t-structures}

We begin by recalling the definitions of the different types of orthogonal collections we will use in this article. The most important are \emph{simple-minded collections}, which model sets of simple objects in the hearts of bounded t-structures in triangulated categories, see \cite{Al-Nofayee,Koenig-Yang}, and \emph{simple-minded systems}, which model sets of nonprojective simple modules in the stable module category of a selfinjective algebra, see \cite{Koenig-Liu, CS15, CSP20}.
In the following we do not require that the collection $\smfull$ is finite unless explicitly stated.

\begin{definition}\label{def:orth}
A collection of objects $\smfull$ in $\sD$ is called \emph{orthogonal} or a \emph{semibrick} if 
\[
\Hom_\sD(u_1,u_2) = 
\begin{cases}
0 & \text{ if } u_1 \not\simeq u_2, \\
\dd_u & \text{ if } u_1 \simeq u_2 (=: u),
\end{cases}
\]
where $\dd_u$ is a division ring.
Let $w \geq 1$ be an integer. An orthogonal collection $\smfull$ is called
\begin{enumerate}
\item (if $w>1$) \emph{$w$-orthogonal} if $\Hom_\sD(u_1 [m], u_2) = 0$ for each $u_1, u_2 \in \smfull$ and each $0< m < w$;
\item a \emph{$w$-simple-minded system} (or \emph{$w$-SMS}) if it is $w$-orthogonal and 
\[
\sD = \extn{\smfull}[w-1] * \cdots * \extn{\smfull}[1] * \extn{\smfull};
\]
\item \emph{$\infty$-orthogonal} if $\Hom_\sD(u_1 [m], u_2) = 0$ for each $u_1, u_2 \in \smfull$ and each $m > 0$;
\item a \emph{simple-minded collection} (or \emph{SMC}) if it is $\infty$-orthogonal and 
\begin{equation} \label{smc-t-structure}
\tag{$\ast$} \sD = \thick{\sD}{\smfull} = \bigcup_{i \geq j} \extn{\smfull}[i] * \cdots * \extn{\smfull}[j].
\end{equation}
\end{enumerate}
\end{definition}

We note that, by \cite[Lemma 2.7]{Dugas}, the extension closure $\extn{\sU}$ of an orthogonal collection $\sU$ is an additive subcategory of $\sD$.

The concepts of $w$-SMS and SMC are closely related to torsion pairs in triangulated categories. 

\begin{definition} \label{def:torsion-pair}
A pair of full additive subcategories $(\sX, \sY)$ of $\sD$ is a \emph{torsion pair} if $\Hom_\sD(\sX,\sY) = 0$ and $\sD = \sX * \sY$. 

If $\sX[1] \subseteq \sX$, then $(\sX,\sY)$ is called a \emph{t-structure} and its \emph{heart} $\sH = \sX \cap \sY[1]$ is an abelian category.
If $\sX[-1] \subseteq \sX$, then $(\sX,\sY)$ is called a \emph{co-t-structure}.

Given $d \in \sD$, the triangle $x \to d \to y \to x[1]$ coming from $\sD =  \sX * \sY$ is called an \emph{approximation triangle} (or \emph{truncation triangle}) for $d$. The triangle is said to be \emph{minimal} if the morphisms $x \to d$ and $d \to y$ are right and left minimal, respectively.

A torsion pair $(\sX,\sY)$ is \emph{bounded} if $\sD = \bigcup_{i \in \bZ} \sX[i] = \bigcup_{i \in \bZ} \sY[i]$. If $(\sX,\sY)$ is a bounded t-structure with heart $\sH$, then, see e.g. \cite[Lem.~3.2]{Bridgeland},
\begin{align*}
\sX & = \susp{\sD}{\sH}  = \bigcup_{i \geq 0} \sH[i] * \cdots * \sH[1] * \sH
\text{ and } \\
\sY & = \cosusp{\sD}{\sH[-1]} = \bigcup_{j \leq -1} \sH[-1] * \sH[-2] * \cdots * \sH[j].
\end{align*}
As a consequence, if $\sK$ is the heart of another bounded t-structure and $\sK \subseteq \sH$ then $\sK = \sH$.

Suppose $(\sW,\sX)$ and $(\sX, \sY)$ are torsion pairs. We say that $(\sW,\sX)$ is \emph{left adjacent} to $(\sX,\sY)$ and $(\sX, \sY)$ is \emph{right adjacent} to $(\sW, \sX)$; cf. \cite[Def.~4.4.1]{Bondarko}.
\end{definition}

For an SMC, $\smfull$, in $\sD$, the formula \eqref{smc-t-structure} for $\thick{\sD}{\smfull}$ is a consequence of the following theorem.
Recall that an abelian category $\sH$ is \emph{length} if it is artinian and noetherian, i.e. each object has a finite composition series whose factors are simple objects of $\sH$.

\begin{theorem}[{\cite[Thm.~4.6]{Schnuerer}}] \label{thm:length-t-structure}
Let $\sD$ be a triangulated category. Then there is a bijection between SMCs in $\sD$ and bounded t-structures whose heart is a length category. 
The map is given by $\smfull \mapsto (\susp{\sD}{\smfull}, \cosusp{\sD}{\smfull[-1]})$; 
the heart is $\extn{\smfull}$.
\end{theorem}

We call a t-structure $(\sX,\sY)$ in $\sD$ \emph{length} if its heart $\sH$ is a length category. In this case, we refer to $\sH$ as a \emph{length heart}. By Theorem~\ref{thm:length-t-structure}, if $(\sX,\sY)$ is a length t-structure, there exists an SMC $\smfull$ such that $(\sX, \sY) = (\susp{\sD}{\smfull},\cosusp{\sD}{\smfull[-1]})$.
If the SMC $\smfull$ is finite then the t-structure and its heart are called \emph{algebraic}.

\subsection{Stability conditions}

We recall the following material on stability conditions from \cite{Bridgeland, Kontsevich-Soibelman}.
The material in this section is only needed in order to be able to understand condition \ref{phase-gap} of Theorem~\ref{thm:length-simple-tilts} and can safely be skipped over in order to follow the rest of the article.

Let $\sH$ be an abelian category and $\lambda \colon K_0(\sH) \to \Lambda$ be a surjective group homomorphism onto a finite rank lattice. 
Let $\Hminus = \{ re^{i \pi \phi} \mid r >0, \phi \in (0,1]\}.$
A \emph{stability function on $\sH$} consists of a group homomorphism $Z \colon \Lambda \to \bC$ such that $Z(h) \coloneqq Z(\lambda(h)) \in \Hminus$ for each $h$ in $\sH$.
If $Z(h) = m e^{i \pi \phi} \in \bC$, then the \emph{phase of $h$} is $\phi(h) = \phi$.
An object $h \in \sH$ is \emph{semistable} if $\phi(h') \leq \phi(h)$ for all nonzero subobjects $h' \into h$, or equivalently, if $\phi(h) \leq \phi(h'')$ for all nonzero quotients $h \onto h''$.
The full subcategory $\cP(\phi)$ of semistable objects of phase $\phi$ is an abelian subcategory of $\sH$.

A stability function $Z$ satisfies the \emph{Harder--Narasimhan} (or \emph{HN}) \emph{property}  if for each $0 \neq h \in \sH$ there exists a filtration
\[
0 = h_0 \into h_1 \into \cdots \into h_{n-1} \into h_n = h
\]
such that $h_i/h_{i-1} \in \cP(\phi_i)$ and $\phi_1 > \phi_2 > \cdots > \phi_n$.
If $\sH$ is length, then any stability function on $\sH$ satisfies the HN property by \cite[Prop.~2.5]{Bridgeland}.

Fix an inner product on $\Lambda_\bR \coloneqq \Lambda \otimes \bR$ and let $\lVert \cdot \rVert$ be the associated norm. 
A stability function $Z$ satisfies the \emph{support property} if there is $C > 0$ with $\lvert Z(h) \rvert \geq C \lVert \lambda(h) \rVert$ for each semistable object $h$ of $\sH$.

\begin{definition}[{\cite[Prop.~5.3]{Bridgeland}}, Stability condition]
A \emph{stability condition} on $\sD$ consists of a pair $(Z, \sH)$ in which $\sH$ is the heart of a bounded t-structure in $\sD$ and $Z$ is a stability function on $\sH$ satisfying the HN and support properties.
\end{definition}

For $\phi \in \bR$, write $\phi = n + \phi_0$ with $\phi_0 \in (0, 1]$ and define $\cP(\phi) \coloneqq \cP(\phi_0)[n]$. The collection of the subcategories $\{ \cP (\phi) \mid \phi \in \bR\}$ is called a \emph{slicing}. 
For an interval $I \subseteq \bR$, define $\cP(I)$ to be the extension closure of $\cP(\phi)$ for $\phi \in I$. We abuse notation by writing $\cP(a,b)$ for $\cP((a,b))$ and so on.
Note that $\cP(0,1] = \sH$.

The support property implies that $\cP(I)$ is a length category  whenever the length of the interval $I$ is strictly less than $1$, see \cite[Lem. 4.5]{Bridgeland08}. We remark that what was referred to as a stability condition in \cite{Bridgeland} is now often referred to as a `pre-stability condition'. The term `stability condition' is usually reserved for pre-stability conditions whose stability functions satisfy the support property because this implies that the associated slicing is locally finite, i.e. each $\cP(\phi)$ is a length category; see \cite[\S 2.1]{Kontsevich-Soibelman}.

\section{Simple tilts to length hearts}

Let $\sH$ be a length heart in $\sD$. In this section we characterise when a simple Happel--Reiten-Smal\o\ tilt of $\sH$ is again length and relate this to the Koenig--Yang mutation formula for simple-minded collections. We first recall Happel--Reiten--Smal\o\ tilting and simple tilts.

\subsection{Happel-Reiten-Smal\o\ tilting and simple tilts}

Let $\sH$ be an abelian category. A torsion pair in $\sH$ is a pair of subcategories $\tt = (\cT, \cF)$ such that $\Hom_\sH(\cT,\cF) = 0$ and 
\[
\sH = \cT * \cF = \{ h \in \sH \mid \text{there is a s.e.s. } 0 \to t \to h \to f \to 0 \text{ with } t\in \cT \text{ and } f \in \cF\}.
\]
The $*$-product in $\sH$ should cause no confusion with the $*$-product in $\sD$ because if $\sH$ is a heart in $\sD$, one has $\Ext^1_\sH(h_1,h_2) = \Hom_\sD(h_1, h_2[1])$ for all objects $h_1$ and $h_2$ of $\sH$. Therefore, the exact structure on $\sH$ is the restriction of the triangulated structure on $\sD$.

\begin{definition}[HRS tilting, {\cite[Prop.~2.1]{HRS}}]\label{def:HRS}
Let $(\sX,\sY)$ be a t-structure in a triangulated category $\sD$ with heart $\sH = \sX \cap \sY[1]$. 
The \emph{right HRS-tilt} of $(\sX, \sY)$ at a torsion pair $\tt = (\cT,\cF)$ is the t-structure $(\sX * \cT[-1],(\cF * \sY)[-1])$ with heart $\rhrs{\tt}{\sH} = \cF * \cT[-1]$.

Note that $\cT =  \rhrs{\tt}{\sH}[1] \cap \sH$ and $\cF = \rhrs{\tt}{\sH} \cap \sH$.
  
The \emph{left HRS-tilt} of $(\sX, \sY)$ at $(\cT,\cF)$ is defined dually and has heart $\lhrs{\tt}{\sH} = \cF[1] * \cT$.
\end{definition}

\begin{figure}
\begin{center}
\begin{tikzpicture}[scale=0.8]

\fill[gray!65] (-6,0) -- (-0.5,0) -- (-0.5,2) -- (-6,2) -- cycle;
\fill[gray!15] (2,0) -- (7.5,0) -- (7.5,2) -- (2,2) -- cycle;

\fill[gray!55] (-0.5,0) -- (0.75,0) -- (0.75,2) -- (-0.5,2) -- cycle;
\fill[gray!40] (0.75,0) -- (2,0) -- (2,2) -- (0.75,2) -- cycle;

\fill[gray!25] (2,0) -- (3.25,0) -- (3.25,2) -- (2,2) -- cycle;

\draw[decoration={brace,mirror},decorate] (3.25,-0.5) --  (7.5,-0.5);
\draw[decoration={brace,mirror},decorate] (-6,-0.5) -- (3.25,-0.5);
\draw[decoration={brace},decorate] (-6,2.5) -- (2,2.5);
\draw[decoration={brace},decorate] (2,2.5) -- (7.5,2.5);

\draw (-6,2) -- (7.5,2);
\draw (-6,0) -- (7.5,0);

\draw[thick] (-5.5,0) -- (-5.5,2);
\draw[thick] (-3,0) -- (-3,2);
\draw[thick] (-0.5,0) -- (-0.5,2);
\draw[thick] (2,0) -- (2,2);
\draw[thick] (4.5,0) -- (4.5,2);
\draw[thick] (7,0) -- (7,2);

\node at (-4.25,1.5) {$\sH[2]$};
\node at (-1.75,1.5) {$\sH[1]$};
\node at (0.15,0.5) {${\scriptstyle \cT}$};
\node at (0.75,1.5) {$\sH$};
\node at (1.4,0.5) {${\scriptstyle \cF}$};
\node at (2.65,0.5) {${\scriptstyle \cT[-1]}$};
\node at (3.25,1.5) {$\sH[-1]$};
\node at (3.85,0.5) {${\scriptstyle \cF[-1]}$};
\node at (5.75,1.5) {$\sH[-2]$};
\node at (-1.375,-1) {$\sX * \cT[-1]$};
\node at (-2,3) {$\sX$};
\node at (4.75,3) {$\sY$};
\node at (5.375,-1) {$(\cF * \sY)[-1]$};
\end{tikzpicture}
\end{center}
\caption{Schematic showing the t-structure $(\sX, \sY)$ and the right HRS tilted t-structure $\big(\sX * \cT[-1],(\cF * \sY)[-1]\big)$ at the torsion pair $(\cT,\cF)$ in the heart $\sH = \sX \cap \sY[1]$.} \label{fig:schematic}
\end{figure}

The next lemma tells us that extension closures of simple objects induce torsion pairs in length hearts. This enables us to define simple tilts.

\begin{lemma}[{\cite[Thm.~2.11]{CSP20} \& \cite[Thm.~3.3]{Dugas}}] \label{lem:simples-are-funct-finite}
Let $\smfull$ be an orthogonal collection of objects in $\sD$ and $\smsub \subseteq \smfull$. Then $\extn{\smsub}$ is functorially finite in $\extn{\smfull}$. 
In particular, if $\sH$ is a length heart in $\sD$ with simple objects $\smfull$ then $(\extn{\smsub}, \sH \cap \smsub\orth)$ and $({}\orth \smsub \cap \sH, \extn{\smsub})$ are torsion pairs in $\sH$.
\end{lemma}

\begin{definition}[Simple tilts]
Suppose $\sH$ is a length heart in $\sD$ and $\smsub$ is a subset of the simple objects of $\sH$. If $\tt = (\cT,\cF) := (\extn{\smsub}, \sH \cap \smsub\orth)$, then we write $\rhrs{\tt}{\sH} = \rhrs{\smsub}{\sH}$ and say that it is the \emph{right simple tilt of $\sH$ at $\smsub$}.
Similarly, if $\tt = (\cT,\cF) := ({}\orth \smsub \cap \sH, \extn{\smsub})$ then we write $\lhrs{\tt}{\sH} = \lhrs{\smsub}{\sH}$ and say that it is the \emph{left simple tilt of $\sH$ at $\smsub$}.
\end{definition}

In the definition above we impose no requirement that $\smsub$ contains only one object, or even finitely many objects.

The following observation is well known; we include a proof for the reader's convenience.

\begin{lemma} \label{lem:shifted-simples}
Let $\sH$ be a length heart in $\sD$ and  $\smsub$ be a subset of its simple objects. Let $\sK = \rhrs{\smsub}{\sH}$ be the right simple tilt of $\sH$ at $\smsub$. Then $s[-1]$ is simple in $\sK$ for each $s \in \smsub$.
\end{lemma}

\begin{proof}
Let $s \in \smsub$ and write $\tt = (\cT,\cF) := (\extn{\smsub}, \sH \cap \smsub\orth)$.
We need to show that $s[-1]$ has no nontrivial subobjects in $\sK$. To that end, suppose it does and consider the short exact sequence in which the first morphism is assumed to be nonzero,
\begin{equation} \label{ses}
0 \to k' \into s[-1] \onto k'' \to 0
\end{equation}
in $\sK$. As $(\cF, \cT[-1])$ is a torsion pair in $\sK$ and torsionfree classes are closed under subobjects, $k' \in \cT[-1]$. Write $k' = t[-1]$ for some $t \in \cT$. Then, as $s$ is simple in $\sH$, the morphism $t \to s$ must be an epimorphism in $\sH$, in particular, there is another short exact sequence
\[
0 \to a \into t \onto s \to 0
\]
in $\sH$. In $\sD$, this is the rotation of the triangle corresponding to \eqref{ses}, so we have a distinguished triangle in $\sD$,
\[
k' \to s[-1] \to k'' = a \to k'[1].
\]
Therefore, $a \in \sH \cap \sK = \cF$ (see Definition~\ref{def:HRS}). On the other hand, $a \in \cT$ because $t, s \in \cT$ and $\cT = \extn{\smsub}$ is a Serre subcategory. Hence $ k'' = a = 0$. That is $k' \into s[-1]$ is an isomorphism and $s[-1]$ is simple in $\sK$.
\end{proof}

\subsection{Mutation of simple-minded collections}

We next define the simple-minded mutation formulas due to Koenig--Yang in \cite{Koenig-Yang} in the case of simple-minded collections and due to Alex Dugas in \cite{Dugas} in the case of simple-minded systems.

\begin{definition}\label{def:mutation}
Let $\smsub$ be a collection of objects in $\sD$ and let $d$ be an object of $\sD$.
\begin{enumerate}
\item The object $d$ \emph{admits a right mutation with respect to $\smsub$} if $d \in \smsub$ or there is a minimal right $\extn{\smsub}$-approximation $\alpha_d \colon s_d \to d[1]$. In this case the \emph{right mutation, $\rsms{\smsub}{d}$, of $d$} is $d$ itself when $ d\in \smsub$ or otherwise given by the triangle
\begin{equation} \label{right-mutation}
\tag{R} \rsms{\smsub}{d}[-1] \too s_d \rightlabel{\alpha_d} d[1] \too \rsms{\smsub}{d}.
\end{equation}
\item The object $d$ \emph{admits a left mutation with respect to $\smsub$} if $d \in \smsub$ or there is a minimal left $\extn{\smsub}$-approximation $\alpha^d \colon d[-1] \to s^d$. In this case the \emph{left mutation, $\lsms{\smsub}{d}$, of $d$} is $d$ itself when $ d\in \smsub$ or otherwise given by the triangle
\begin{equation} \label{left-mutation}
\tag{L} \lsms{\smsub}{d} \too d[-1] \rightlabel{\alpha^d} s^d \too \lsms{\smsub}{d}[1].
\end{equation}
\end{enumerate}
Given a collection $\smfull$ of objects in $\sD$, we set $\rsms{\smsub}{\smfull} \coloneqq \{\rsms{\smsub}{u} \mid u \in \smfull\}$. The notation $\lsms{\smsub}{\smfull}$ is defined similarly. 
\end{definition}

Note that the definition of mutation of simple-minded collections introduced in \cite[Def. 7.5]{Koenig-Yang} differs from the definition above by a shift. More precisely, if $\sU$ is an SMC and $\sS \subseteq \sU$, then the right mutation of $\sU$ at $\sS$ as defined in \cite{Koenig-Yang} coincides with $\rsms{\smsub}{\sU}[-1]$, as defined above.
In particular, if $\sH = \extn{\smfull}$ then we will observe in \ref{SMC} $\implies$ \ref{length} in Theorem~\ref{thm:length-simple-tilts} below that $\rhrs{\smsub}{\sH} = \extn{\rsms{\smsub}{\smfull}[-1]}$.

The following property of simple-minded right approximations is useful; there is a dual result for left approximations.
When $\sD$ is Hom-finite and Krull--Schmidt, the right approximations exist if and only if minimal right approximation exist, see e.g. \cite{Auslander-Smalo}. 

\begin{lemma}[{\cite[Lem.~4.6]{Dugas}}] \label{lem:Dugas}
Let $\smsub \subseteq \sD$ be an orthogonal collection, and suppose $d \in \sD$ admits a right $\extn{\smsub}$-approximation. Consider the minimal $(\extn{\smsub}, \smsub\orth)$-triangle for $d$:
\begin{equation} \label{mutation-triangle}
s_d \rightlabel{f} d \too z_d \too s_d[1].
\end{equation}
Then the map $\Hom(\smsub,f) \colon \Hom(\smsub, s_d) \to \Hom(\smsub, d)$ is an isomorphism.
\end{lemma}

\subsection{Characterisation of length simple tilts}

We can now formulate the characterisation of when simple tilts are length.

\begin{theorem} \label{thm:length-simple-tilts}
Let $\sH$ be a length heart in $\sD$. Suppose $\smfull = \simples{\sH}$ is the set of simple objects of $\sH$ and let $\smsub \subseteq \smfull$. 
Let $\tt = (\cT,\cF) := (\extn{\smsub}, \smsub\orth\cap \sH)$ and write $\sK = \rhrs{\smsub}{\sH} = \cF * \cT[-1]$ for the right simple tilt of $\sH$ at $\smsub$. 
Then the following conditions are equivalent.
\begin{enumerate}[label=(\roman*)]
\item $\sK$ is a length heart. \label{length}
\item Each object in $\sK[1]$ admits a right $\extn{\smsub}$-approximation.  \label{K-approx}
\item Each object in $\cF[1]$ admits a right $\extn{\smsub}$-approximation.  \label{F-approx}
\item Each object in $(\smfull \setminus \smsub)[1]$ admits a right $\extn{\smsub}$-approximation. \label{T-approx}
\item The right mutation $\rsms{\smsub}{\smfull}$ exists and is a simple-minded collection. \label{SMC} 
\end{enumerate}
Suppose further that $\sH$ is length with finitely many isoclasses of simple objects. Then the conditions above are equivalent to the following condition.
\begin{enumerate}[resume, label=(\roman*)]
\item There exists a stability condition $\sigma = (Z,\cP)$ with $\cP(0,1] = \sH$, $\cP(1)=\extn{\smsub}$ and $\cP(0,\phi] = \smsub\orth\cap \sH$ for some $0<\phi<1$. \label{phase-gap}
\end{enumerate}
\end{theorem}

\begin{proof}
\ref{length} $\implies$ \ref{K-approx}. Suppose $\sK$ is length. By Lemma~\ref{lem:shifted-simples}, $\smsub[-1]$ is a subset of the simple objects of $\sK$. Therefore, $\extn{\smsub[-1]}$ is functorially finite in $\sK$ by Lemma~\ref{lem:simples-are-funct-finite} as $\sK$ is length. Hence, $\extn{\smsub}$ is functorially finite in $\sK[1]$. In particular, each object of $\sK[1]$ admits a right $\extn{\smsub}$-approximation.

\ref{K-approx} $\implies$ \ref{F-approx} $\implies$ \ref{T-approx}. We have $(\smfull\setminus \smsub) \subseteq \cF \subseteq \sK$.

\ref{T-approx} $\implies$ \ref{SMC}. This is the content of \cite[Prop.~7.6(c)]{Koenig-Yang}. The first condition in {\it loc.~cit.} is our assumption \ref{T-approx}. The second condition follows from Lemma~\ref{lem:Dugas}. We observe that Koenig and Yang's third condition is not required: the injectivity of the map $\Hom(\smsub,f[1]) \colon \Hom(\smsub,  s_d[1]) \to \Hom(\smsub, d[1])$ induced by \eqref{mutation-triangle} is required to get the vanishing of $\Hom(\smsub,z_d)$, which follows from the triangulated Wakamatsu lemma, see e.g. \cite[Lem. 2.1]{J09}.

For the convenience of the reader, we reproduce Koenig and Yang's proof below and remark on where our hypotheses are used in relation to Koenig and Yang's hypotheses. 
Koenig and Yang's original argument assumes that $\smsub = \{ s \}$. We note that this is not required for the argument to work; indeed, $\smsub$ may even be infinite.

Let $u_1, u_2 \in \smfull \setminus \smsub$. By the assumption that each object of $(\smfull \setminus \smsub)[1]$ admits a right $\extn{\smsub}$-approximation and every right approximation contains a minimal right approximation as a summand, there are mutation triangles
\[
s_{u_i} \rightlabel{\alpha_i} u_i[1] \too \rsms{\smsub}{u_i} \too s_{u_i}[1] \quad \text{for } i =1,2.
\]
These show that $\rsms{\smsub}{u_i} \in \smfull[1] * \extn{\smsub}[1]$ from which we obtain immediately that
\begin{align*}
\Hom(\rsms{\smsub}{u_1}[>0], \rsms{\smsub}{u_2}) & = 0, \\
\Hom(\rsms{\smsub}{u_i}[\geq 0], \smsub)         & = 0, \\
\Hom(\smsub[>1], \rsms{\smsub}{u_i})             & = 0.
\end{align*}
It remains to check that $\Hom(\smsub,\rsms{\smsub}{u_i}) = 0$,  $\Hom(\smsub[1],\rsms{\smsub}{u_i}) = 0$, $\Hom(\rsms{\smsub}{u_1}, \rsms{\smsub}{u_2}) \simeq \Hom(u_1, u_2)$ and that $\rsms{\smsub}{\smfull}$ generates $\sD$.

Write $u = u_1$ and $\alpha = \alpha_1$.
The first claim follows from the Wakamatsu lemma (this is where Koenig and Yang require that $\Hom(\smsub,\alpha[1])$ is injective).
The second claim follows from applying $\Hom(\smsub,-)$ to the mutation triangle for $u$ to get,
\[
(\smsub, u) \too (\smsub, \rsms{\smsub}{u}[-1]) \too (\smsub, s_u) \rightlabel{(\smsub,\alpha)} (\smsub,u[1]).
\]
Here, Koenig and Yang use $\Hom(\smsub, u) = 0$ and the injectivity of $\Hom(\smsub,\alpha)$ to get that $\Hom(\smsub[1],\rsms{\smsub}{u}) = 0$. For us, Lemma~\ref{lem:Dugas}, tells us that $\Hom(\smsub,\alpha)$ is injective.

For the third claim, we apply $\Hom(u_1,-)$ to the mutation triangle for $u_2$,
\[
(u_1,s_{u_2}[-1]) \to (u_1,u_2) \to (u_1,\rsms{\smsub}{u_2}[-1]) \to (u_1, s_{u_2}),
\]
where the vanishing of the outer terms forces $\Hom(u_1,u_2) \simeq \Hom(u_1,\rsms{\smsub}{u_2}[-1])$.
Next, applying $\Hom(-,\rsms{\smsub}{u_2}[-1])$ to the rotated mutation triangle for $u_1$ gives,
\[
(s_{u_1},\rsms{\smsub}{u_2}[-1]) \to (\rsms{\smsub}{u_1}[-1],\rsms{\smsub}{u_2}[-1]) \to (u_1, \rsms{\smsub}{u_2}[-1]) \to (s_{u_1}[-1],\rsms{\smsub}{u_2}[-1]),
\]
where the vanishing of the outer terms again forces
\begin{align*}
\Hom(\rsms{\smsub}{u_1},\rsms{\smsub}{u_2}) & \simeq \Hom(\rsms{\smsub}{u_1}[-1],\rsms{\smsub}{u_2}[-1]) \\
                                                                           & \simeq \Hom(u_1,\rsms{\smsub}{u_2}[-1]) \\
                                                                           & \simeq \Hom(u_1,u_2).
\end{align*}

Finally, one observes immediately that the mutation triangles show generation. Hence, $\rsms{\smsub}{\smfull}$ is an SMC.

\ref{SMC} $\implies$ \ref{length}. It is sufficient to show that $\extn{\rsms{\smsub}{\smfull}[-1]} \subseteq \sK$. As both $\extn{\rsms{\smsub}{\smfull}}$ and $\sK$ are hearts of bounded t-structures in $\sD$, it follows that $\extn{\rsms{\smsub}{\smfull}[-1]} = \sK$, see Definition~\ref{def:torsion-pair}. Therefore, $\sK$ is a length heart whose simple objects are exactly those objects in the SMC, $\rsms{\smsub}{\smfull}[-1]$.
Note that $\sK = \cF * \cT[-1] = \extn{\cF \cup \cT[-1]}$ since $\cT[-1] \subseteq \sK$, $\cF \subseteq \sK$ and $\sK$ is the heart of a bounded t-structure and therefore closed under extensions.
Clearly, $\smsub[-1] \subseteq \cT[-1] = \extn{\smsub}[-1]$. 
For $u \in \smfull \setminus \smsub$, applying $\Hom(\smsub,-)$ to the mutation triangle \eqref{right-mutation} and observing that $\Hom(\smsub[1], \rsms{\smsub}{u}) = 0$ by Lemma~\ref{lem:Dugas}, gives $\rsms{\smsub}{u} \in \sH[1] \cap (\smsub[1])\orth = \cF[1]$. That is, $\rsms{\smsub}{u}[-1] \in \cF$ and $\rsms{\smsub}{\smfull}[-1] \subseteq \cF \cup \cT[-1]$, whence  $\extn{\rsms{\smsub}{\smfull}[-1]} \subseteq \sK$, as required.

Now suppose, in addition, that $\sH$ has finitely many isoclasses of simple objects. In this case, we show that \ref{length} $\iff$ \ref{phase-gap}. 

\ref{length} $\implies$ \ref{phase-gap}. Since $\sH$ is length, stability conditions $\sigma=(Z,\cP)$ with $\cP(0,1]=\sH$ correspond bijectively to stability functions $Z$  mapping the simple objects of $\sH$ into $\Hminus$. Define a stability function $Z$ by setting $Z(s)=-1$ for $s\in \smsub$ and $Z(u)=i$ for $u\in \smfull\setminus \smsub$. Then $\cP(1)=\extn{\smsub}$ and $\cP(0,1) = \cP(1)\orth = \smsub\orth \cap \sH$. Hence, 
\[
\sK=(\smsub\orth\cap \sH) * \extn{\smsub}[-1]= \cP(0,1)*\cP(0) = \cP[0,1).
\]
Since $\sK$ is length and the number of isoclasses of simple objects of a length heart is equal to the rank of the Grothendieck group, $\sK$ has finitely many isoclasses of simple objects. These generate $\sK$ by extensions. It follows that $\sK = \cP[0,\phi]$ for some $0<\phi <1$, namely for $\phi$ the maximal phase of any HN factor of a simple object in $\sK$ with respect to the slicing $\cP$. Thus, $\cP(\phi,1)=\{0\}$ and $\smsub\orth\cap \sH = \cP(0,1) = \cP(0,\phi]$ as required.

\ref{phase-gap} $\implies$ \ref{length}. Let $\sigma=(Z,\cP)$ be a stability condition with heart $\cP(0,1]=\sH$, $P(1)=\extn{\smsub}$ and $\cP(0,\phi]=\smsub\orth\cap\sH$ for some $0< \phi<1$. Then 
\[
\sK = (\smsub\orth\cap\sH) * \extn{\smsub}[-1] = \cP(0,\phi] * \cP(0) = \cP[0,\phi].
\]
Since $[0,\phi]$ is an interval whose length is strictly smaller than $1$, the heart $\sK$ is length by \cite[Lem. 4.5]{Bridgeland08}.
\end{proof}

\begin{remark}
There is an evident result dual to Theorem~\ref{thm:length-simple-tilts} for left simple tilts/left simple-minded mutation.
\end{remark}

Finally, we examine the dual of the counterexample of Koenig and Yang \cite[p.~428]{Koenig-Yang} 
to illustrate that simple tilts of algebraic hearts need not be algebraic 
using the language of Theorem~\ref{thm:length-simple-tilts}. 

\begin{example} \label{ex:koenig-yang}
Consider the quiver $Q$ below.
\[
\begin{tikzcd}
\arrow[loop, out=140, in=220,looseness=5] 1 \ar[r] & 2
\end{tikzcd}
\]
Let $\sD = \Db(\nil{Q})$ be the bounded derived category of nilpotent $\kk$-representations of $Q$, i.e. those whose composition factors belong to the set of the simple representations corresponding to the vertices of $Q$. 
The heart $\sH = \nil{Q}$ is length with two simple objects, $s_1$ and $s_2$, corresponding to each of the vertices. 
Consider the right simple tilt $\sK = \rhrs{s_1}{\sH} = (s_1\orth \cap \sH) * \extn{s_1} [-1]$.
We use Theorem~\ref{thm:length-simple-tilts} to show that $\sK$ is not length, first by using stability conditions and second by using approximation theory.

For the stability condition approach, consider the uniserial nilpotent object $m_n$ for $n \geq 1$ in which $s_1$ occurs as a composition factor $n$ times and $s_2$ occurs as a composition factor exactly once:
\[
m_n = 
\begin{matrix} s_1\\[-2pt] \vdots \\[-2pt] s_1 \\[-2pt] s_2 \end{matrix}.
\] 
As $\sH$ is length, stability conditions $\sigma = (Z,\cP)$ such that $\cP(1) = \extn{s_1}$ correspond bijectively with stability functions $Z$ such that $Z(s_1) \in \bR_{<0}$. 
Without loss of generality, we may fix $Z(s_1) = -1$.
In order to satisfy $\cP(0,\phi] = s_1\orth \cap \sH$, we require that $0 < \phi(s_2) < 1$. 
For any choice of $Z(s_2)$ such that $0 < \phi(s_2) < 1$, we have that $m_n$ is $\sigma$-semistable and $\lim_{n \to \infty} \phi(m_n) = 1$. As $m_n$ has simple socle $s_2$, we have $m_n \in s_1\orth \cap \sH$.
Hence, there is no stability condition $\sigma = (Z,\cP)$ such that $\cP(1) = \extn{s_1}$ and $\cP(0,\phi] = s_1\orth \cap \sH$ and by Theorem~\ref{thm:length-simple-tilts}\ref{phase-gap}, 
$\sK$ is not length.

For the approximation theory approach,
suppose $s_2[1]$ admits a right $\extn{s_1}$-approximation, $f \colon x \to s_2[1]$. Each indecomposable object 
\[
x_n = 
\begin{matrix} s_1\\[-2pt] \vdots \\[-2pt] s_1 \end{matrix} \in \extn{s_1}
\] 
admits a nonzero homomorphism $g \colon x_n \to s_2[1]$, whose cocone is the object $m_n$ above. This morphism must factor through $f$, i.e.~$g= fh$ for some $h\colon x_n \to x$, giving rise to the commutative diagram coming from the octahedral axiom below, where $k = \ker h$, $c = \coker h$ and the cone of $h$ has the given form because $\nil{Q}$ is hereditary.
\[
\begin{tikzcd}
s_2 \ar{r} \ar[equals]{d} & m_n \ar{r} \ar{d}             & x_n \ar{d}{h} \\
s_2 \ar{r}                & u \ar{r} \ar{d}               & x \ar{d} \\
                          & c \oplus k[1] \ar[equals]{r}  & c \oplus k[1]
\end{tikzcd}
\]
In particular, if $k \neq 0$, then the middle column says that $s_1$ is a simple subobject of $m_n$ as $k \in \extn{s_1}$. This contradicts the fact that $m_n$ is uniserial with simple socle $s_2$. Hence $k = 0$, meaning the factoring map $h\colon x_n \to x$ must be injective. But since the length of $x_n$ is arbitary, choosing $n$ larger than the length of $x$ gives a contradiction. Hence, there is no such right $\extn{s_1}$-approximation of $s_2[1]$.
\end{example}

\section{Infinitely iterable simple-minded mutation} \label{sec:iteration}

Let $\smfull$ be an SMC and $\smsub \subseteq \smfull$. In the previous section, we obtained conditions equivalent to the right mutation $\rsms{\smsub}{\smfull}[-1]$ also being an SMC. 
However, this process may not be iterable; see Example~\ref{ex:koenig-yang-2} below.
The main result of this section states that iterated mutations of nice enough simple-minded systems/collections are simple-minded systems/collections. Here, `nice enough' means that the SMS/SMC Setups \ref{SMC-setup} and \ref{SMS-setup} are satisfied. 

Let $\smfull$ be an SMC or a $w$-SMS for $w \geq 1$, and suppose $\smsub \subseteq \smfull$. Define $\rsmsn{1}{\smsub}{\smfull} = \rsms{\smsub}{\smfull}$ and for $n > 1$, define
\[
\rsmsn{n}{\smsub}{\smfull} = \rsms{\smsub}{\rsmsn{n-1}{\smsub}{\smfull}},
\]
when they make sense.
One defines $\lsmsn{n}{\smsub}{\smfull}$ for $n \geq 1$ similarly.

\begin{theorem}\label{thm:mut-SM-is-SM}
Let $\sD$ be a Hom-finite, Krull--Schmidt triangulated category.
\begin{enumerate}
\item Suppose $\smfull$ is an SMC and $\smsub \subseteq \smfull$. Then  $\rsmsn{n}{\smsub}{\smfull}$ and $\lsmsn{n}{\smsub}{\smfull}$ are defined and are again SMCs for each $n \geq 1$ if, and only if, $\extn{\smsub}$ is covariantly finite in ${}\orth (\smsub[{<}0])$ and contravariantly finite in $(\smsub[{>}0])\orth$;
\item Suppose $\smfull$ is a $w$-SMS for  an integer $w \geq 1$ and $\smsub \subseteq \smfull$. Then $\rsmsn{n}{\smsub}{\smfull}$ and $\lsmsn{n}{\smsub}{\smfull}$ are defined and are again $w$-SMSs for each $n \geq 1$ if $\extn{\smsub}$ is functorially finite and $\smsub$ is invariant under the shift of the Serre functor, $\SSS[w]$. Moreover, if $\sD$ is $(-w)$-Calabi--Yau, the converse also holds.
\end{enumerate}
\end{theorem}

\begin{remark}
Theorem~\ref{thm:mut-SM-is-SM} provides a unified proof of \cite[Thm.~4.2]{Dugas} for $1$-SMSs and \cite[Thm.~5.3]{J22} for $w$-SMSs, with $w \geq 2$. The latter was approached via an analogue of HRS tilting theory.
\end{remark}

\begin{example} \label{ex:koenig-yang-2}
We adapt Example~\ref{ex:koenig-yang} to give an example of an SMC $\smfull$ with $\smsub \subseteq \smfull$ for which $\rsms{\smsub}{\smfull}$ is again an SMC but $\rsmsn{2}{\smsub}{\smfull}$ is not defined.
Keeping the notation of Example~\ref{ex:koenig-yang}, we set $\sL = \lhrs{s_1}{\sH}[-1] = \extn{s_1} * ({}\orth s_1 \cap \sH)[-1]$ to be the inverse shift of the simple left HRS tilt of $\sH$ at the torsion pair $\tt = ({}\orth s_1 \cap \sH, \extn{s_1})$.

As $\Hom(s_2[-1],s_1) \simeq \Ext^1(s_2,s_1) = 0$, $s_2[-1]$ admits a (trivial) left $\extn{s_1}$-approximation and $\sL$ is algebraic by the dual of Theorem~\ref{thm:length-simple-tilts}.
Applying the mutation formula in Definition~\ref{def:mutation}\eqref{left-mutation}, we find that the simple objects of $\sL$ are $s_1$ and $s_2[-1]$. Further, $\Ext^1(s_2[-1],s_1) = 0$ as $\sH$ is hereditary and $\Ext^1(s_1, s_2[-1]) = 0$ by Schur's lemma, which imply that $\sL = \add{s_2[-1], \extn{s_1}}$. 
Noting left and right tilting are inverse operations, we have $\rhrs{s_1}{\sL}[1] = \sH$.
By Example~\ref{ex:koenig-yang}, we know that $\rhrs{s_1}{\sH}[1] = \rhrs{s_1}{\rhrs{s_1}{\sL}[1]}[1]$ is not length. Therefore,  $\rsmsn{2}{s_1}{\{s_1,s_2[-1]\}}$ is not defined.

Let us show explicitly that the condition equivalent to infinite iterability in Theorem~\ref{thm:mut-SM-is-SM} is also broken.
In Example~\ref{ex:koenig-yang}, we showed that $s_2[1]$ does not admit a right $\extn{s_1}$-approximation. As $\sU = \{s_1, s_2\}$ is an SMC, we have $\Hom(s_1[>0], s_2[1])=0$. In particular, $s_2[1] \in (s_1[>0])\orth$ is an object not admitting a right $\extn{s_1}$-approximaton. As $\extn{s_1} \subseteq (s_1[>0])\orth$, it follows that $\extn{s_1}$ is not contravariantly finite in $(s_1[>0])\orth$.

Finally, setting $\sL_n = \extn{\lsmsn{n}{s_1}{\{s_1,s_2\}}}$ and noting the fact that $\Hom(s_2[-n],s_1)=0$ for each $n > 0$ means that the dual of Theorem~\ref{thm:length-simple-tilts} applies and $\sL_n$ is length with simple objects $\lsmsn{n}{s_1}{\{s_1,s_2\}} =  \{s_1, s_2[-n]\}$ for each $n$. 
Then $\rsmsn{i}{s_1}{\{s_1,s_2[-n]\}}$ is defined and an SMC for each $1 \leq i \leq n$, but $\rsmsn{n+1}{s_1}{\{s_1,s_2[-n]\}}$ is not defined, showing that it is possible to be able to iteratively perform an arbitrary finite number of simple-minded mutations but not infinitely many simple-minded mutations. 
\end{example}

\subsection{Simple minded reduction}

As described in the introduction, the philosophy behind mutation is `reduce, shift, lift'. In order to describe the reduction step in this procedure, we recall the setups required for SMC reduction in \cite{Jin} and SMS reduction \cite{CSP20}.

\begin{setup}[SMC Setup] \label{SMC-setup}
Let $\smsub$ be an $\infty$-orthogonal collection of objects in $\sD$ and $\sZ$ a subcategory of $\sD$ satisfying the following conditions:
\begin{enumerate}
\item $\extn{\smsub}$ is covariantly finite in ${}\orth (\smsub[{<}0])$ and contravariantly finite in $(\smsub[{>}0])\orth$; \label{smc1}
\item for $d \in \sD$, we have $\Hom_\sD (d, \smsub[{\ll} 0]) = 0$ and  $\Hom_\sD (\smsub[{\gg} 0], d) = 0$; and, \label{smc2}
\item $\sZ = {}\orth (\smsub[\leq 0]) \cap (\smsub[\geq 0])\orth$. \label{smc3}
\end{enumerate}
\end{setup}

If $\smfull$ is an SMC and $\smsub \subseteq \smfull$ satisfies Setup~\ref{SMC-setup}, then Theorem~\ref{thm:length-simple-tilts}\ref{T-approx} holds and mutation is defined.

\begin{remark} \label{rem:smc2-is-necessary}
Condition \eqref{smc2} above is necessary for $\smsub$ to be a subset of an SMC.
Indeed, suppose $\smfull$ is an SMC and let $d$ be an object of $\sD$. Then there exist integers $i \geq j$ such that $d \in \extn{\smfull}[i] * \extn{\smfull}[i-1] * \cdots * \extn{\smfull}[j]$. Then,  we have $\Hom_\sD(\smfull[>i],d) = 0$ and $\Hom_\sD(d, \smfull[<j]) = 0$ as $\Hom_{\sD}(\smfull[>0],\smfull) = 0$.
\end{remark}

For $w$-SMSs, we recall the following from \cite{CSP20}.
Let $w \geq 1$ and $\sX$ be a subcategory of $\sD$. We denote by $\sX\orthw$ the following right perpendicular category
\[
\sX\orthw \coloneqq \{d \in \sD \mid \Hom (\sX[i], d) = 0 \text{ for } i= 0, \ldots, w\}. 
\]
The left perpendicular category ${}\orthw \sX$ is defined dually. Assume $\sD$ has a Serre functor $\SSS$. 
\begin{setup}[SMS Setup] \label{SMS-setup}
Let $w \geq 1$. Let $\smsub$ be a $w$-orthogonal collection and $\sZ$ be a subcategory of $\sD$ satisfying the following conditions:
\begin{enumerate}
\item  $\extn{\smsub}$ is functorially finite and $\smsub$ is invariant under $\SSS[w]$; and,
\item $\sZ = \smsub\orthw$. 
\end{enumerate}
\end{setup}

If $\smfull$ is a $w$-SMS and $\smsub \subseteq \smfull$ then \cite[Cor.~2.9]{CSP20} implies $\smsub$ satisfies the functorial finiteness condition in Setup~\ref{SMS-setup} and $\rsms{\smsub}{u}$ and $\lsms{\smsub}{u}$ are well defined for each $u \in \smfull \setminus \smsub$.

The key technical tool used to prove Theorem~\ref{thm:mut-SM-is-SM} is simple-minded reduction, which we recall from \cite[Appendix A]{CSPP} and \cite[Thm.~3.1]{Jin} for SMCs and \cite[Thms.~A \& B]{CSP20} for $w$-SMSs.

\begin{theorem}[Simple-minded reduction] \label{thm:sm-reduction}
Let $\smsub$ and $\sZ$ be as in the SMC Setup~\ref{SMC-setup} (SMS Setup~\ref{SMS-setup}, respectively).
Then $\sZ$ is a triangulated category with shift functor $\shift{1} \colon \sZ \to \sZ$ defined by taking the cone of a minimal right $\extn{\smsub}$-approximation, 
\[
s_z \to z[1] \to z\shift{1} \to s_z[1].
\]
Moreover, there is a bijection
\begin{align*}
\{ \text{SMCs ($w$-SMSs, resp.) in $\sD$ containing $\smsub$} \} & \bij  \{ \text{SMCs ($w$-SMSs, resp.) in $\sZ$} \}. \\
\smfull                                                                                               & \longmapsto  \smfull \setminus \smsub 
\end{align*}
\end{theorem}

Note that the shift functor $\shift{1}$ is defined on morphisms in the obvious way, and its quasi-inverse $\shift{-1}$ by the dual construction, see \cite[Lem.~3.6]{CSP20}. We refer to \cite[Thms.~4.1 \& 5.1]{CSP20} and \cite[Appendix A]{CSPP} for an explicit description of the triangulated structure on $\sZ$.

\begin{remark} \label{rem:sm-reduction}
We highlight the following features of Theorem~\ref{thm:sm-reduction}.
\begin{enumerate}[label=(\roman*)]
\item The triangles used to define $\shift{1}$ and $\shift{-1}$ are `generalised' right and left $\smsub$-mutation triangles. Here `generalised' means we mutate any object of $\sZ$ with respect to $\smsub$ not just objects of $\smfull \setminus \smsub$ for an SMC or $w$-SMS, $\smfull$.
\item In the SMC case there are one-sided triangulated structures. In particular, examining the proof of \cite[Thm.~A.2]{CSPP} shows that replacing  
\eqref{smc1} and \eqref{smc2} in Setup~\ref{SMC-setup} with
\begin{enumerate}
\item[$(1')$] $\extn{\smsub}$ is contravariantly finite in $(\smsub[>0])\orth$; and,
\item[$(2')$] for $d \in \sD$, we have $\Hom_\sD(\smsub[\gg 0],d) =0$,
\end{enumerate}
and leaving \eqref{smc3} the same, gives $\sZ$ the structure of a right triangulated category (see \cite{BM94} for the definition) with shift endofunctor $\shift{1}$. \label{right-triangulated}
\item Dually, taking the other obvious replacements of \eqref{smc1} and \eqref{smc2}, gives $\sZ$ the structure of a left triangulated category with loop endofunctor $\shift{-1}$. \label{left-triangulated}
\item In the case of SMC Setup~\ref{SMC-setup}, the inclusion functor composed with the quotient functor induces a triangle equivalence $\sZ \simeq \sD/\thick{\sD}{\smsub}$ by \cite[Thm.~3.1]{Jin}.
\end{enumerate}
\end{remark}

\subsection{Mutation pairs}

We obtain a proof of Theorem~\ref{thm:mut-SM-is-SM} as a consequence of a statement on simple-minded mutation pairs based on an idea in \cite{IYo}.
We recall the definition from \cite[Def.~3.2]{CSP20}.

\begin{definition}
Let $\smsub$ be a collection of objects of $\sD$ considered as a full subcategory. A pair $(\sA,\sB)$ of full subcategories of $\sD$ is called an \emph{$\smsub$-mutation pair} if 
\[
\sA = {}\orth \smsub\orth \cap {}\orth(\smsub[-1]) \cap (\extn{\smsub} * \sB)[-1]
\quad \text{and} \quad
\sB = {}\orth \smsub\orth \cap (\smsub[1])\orth \cap (\sA * \extn{\smsub})[1],
\]
where ${}\orth \smsub\orth = {}\orth \smsub \cap \smsub\orth$.
\end{definition}

We require the following addendum to Lemma~\ref{lem:Dugas}, which is a generalisation of \cite[Lem.~4.7]{Dugas} (see also \cite[Lem.~2.6]{CSP20}). The same proof carries over in this setting.

\begin{lemma}\label{lem:Dugas-addendum}
Assume the notation and setup of Lemma~\ref{lem:Dugas}.
Suppose further that $w \geq 1$ and $\smsub \subset \sD$ is invariant under $\SSS[w]$. 
\begin{enumerate}[resume]
\item The map $\Hom(g[w-1], \smsub) \colon \Hom(z_d[w-1], \smsub) \to \Hom(d[w-1], \smsub)$ is a monomorphism.
\item If $d \in {}\orth (\smsub[1-w])$ then $z_d \in {}\orth (\smsub[1-w])$.
\end{enumerate}
\end{lemma}

We also require the following small generalisation of \cite[Lem.~3.6]{CSP20} in order to apply it in the SMC context. Again the same proof carries over.

\begin{lemma} \label{lem:pair-implies-mutation}
Let $\smsub$ be an orthogonal collection of objects in $\sD$ satisfying $\smsub$ is invariant under $\SSS[1]$ or $\Hom(\smsub[1], \smsub) = 0$.
Let $(\sA,\sB)$ be an $\smsub$-mutation pair. Then the following hold. 
\begin{enumerate}
\item If each object of $\sA[1]$ admits a right $\extn{\smsub}$-approximation then $\sB = \rsms{\smsub}{\sA}$. \label{right-approx-exists}
\item If each object of $\sB[-1]$ admits a left $\extn{\smsub}$-approximation then $\sA = \lsms{\smsub}{\sB}$. \label{left-approx-exists}
\end{enumerate} 
\end{lemma}

\begin{remark} \label{rem:pair-implies-mutation}
If $\smsub$ is a $w$-orthogonal collection satisfying SMS Setup~\ref{SMS-setup} then the hypotheses in Lemma~\ref{lem:pair-implies-mutation} are satisfied for any $\smsub$-mutation pair.
Now suppose $\smsub$ is an $\infty$-orthogonal collection satisfying the SMC Setup~\ref{SMC-setup} and $(\sA,\sB)$ is an $\smsub$-mutation pair. 
Then $\Hom(\smsub[1],\smsub) = 0$. Additionally, if $\sA \subset (\smsub[\geq 0])\orth$ then the hypothesis of Lemma~\ref{lem:pair-implies-mutation}\eqref{right-approx-exists} holds.
Similarly, if  $\sB \subset {}\orth(\smsub[\leq 0])$ then the hypothesis of Lemma~\ref{lem:pair-implies-mutation}\eqref{left-approx-exists} holds.
\end{remark}

The following proposition establishes a relationship between $\smsub$-mutation pairs and simple-minded mutation of $\smfull$ at $\smsub$.

\begin{proposition} \label{lem:mutation-pair-complement}
Suppose $\smfull$ is an SMC ($w$-SMS, resp.) and $\smsub \subseteq \smfull$ satisfies the SMC Setup~\ref{SMC-setup} (SMS Setup~\ref{SMS-setup}, resp.). 
Then $(\smfull \setminus \smsub, \rsms{\smsub}{\smfull\setminus \smsub})$ and $(\lsms{\smsub}{\smfull \setminus \smsub}, \smfull \setminus \smsub)$ are $\smsub$-mutation pairs.  
\end{proposition}

\begin{proof}
Suppose $\smfull$ is an SMC ($w$-SMS, resp.) in $\sD$ satisfying the SMC Setup~\ref{SMC-setup} (SMS Setup~\ref{SMS-setup}, resp.). Let $\sA \coloneqq \smfull \setminus \smsub$ and $\sB \coloneqq \rsms{\smsub}{\sA}$. We first show that 
\[
\sA = {}\orth \smsub\orth \cap {}\orth(\smsub[-1]) \cap (\extn{\smsub} * \sB)[-1].
\] 
The inclusion  $\sA \subseteq {}\orth \smsub\orth \cap {}\orth(\smsub[-1]) \cap (\extn{\smsub} * \sB)[-1]$ is clear since
\[
\sB  = \rsms{\smsub}{\sA}  = \left\{b \ \bigg\vert \
\parbox{11cm}{there exists $a \in \sA$ and a triangle $s_a \rightlabel{\alpha} a[1] \too b \too s_a[1]$ with $\alpha$ a minimal right $\extn{\smsub}$-approximation}\right\},
\]
from which it follows that $\sA \subseteq \extn{\smsub}[-1] * \sB[-1]$. The inclusion $\sA \subseteq  {}\orth \smsub\orth \cap {}\orth(\smsub[-1])$ is immediate in the case where $\smfull$ is an SMC or a $w$-SMS with $w \geq 2$. In the case where $\smfull$ is a $1$-SMS, the fact that $\sA \subseteq {}\orth(\smsub[-1])$ follows from the fact that $\smsub$ is invariant under $\SSS[1]$. 

For the inclusion ${}\orth \smsub\orth \cap {}\orth(\smsub[-1]) \cap (\extn{\smsub} * \sB)[-1] \subseteq \sA$, take $d \in {}\orth \smsub\orth \cap {}\orth(\smsub[-1]) \cap (\extn{\smsub} * \sB)[-1]$. This sits in the commutative diagram constructed from the octahedral axiom, 
\[
\begin{tikzcd}
                                     & a[1] \ar[equals]{r} \ar{d} & a[1] \ar{d}{0} \\
d[1] \ar[equals]{d} \ar{r} & b  \ar{r} \ar{d}                 & s_1[1] \ar{d} \\
d[1] \ar{r}[swap]{0}         & s_2[1] \ar{r}                    & c[1] 
\end{tikzcd}
\]
in which the top horizontal triangle is given since $d \in  (\extn{\smsub} * \sB)[-1]$ and the left-hand vertical triangle is given by $\sB = \rsms{\smsub}{\sA}$. The morphisms marked $0$ are zero because $d \in {}\orth \smsub$ and $\sA \subseteq {}\orth \smsub$.
It follows that $c \simeq s_1 \oplus a[1] \simeq s_2 \oplus d[1]$.  As $\Hom_{\sD}(d[1],s_1) = 0$, $d[1]$ must be a summand of $a[1]$. Similarly, as $\Hom_{\sD}(a[1],s_2) = 0$ since $\sA \subseteq {}\orth(\smsub[-1])$ by the first part of the proof, we have that $a[1]$ is a summand of $d[1]$. Hence, $d \simeq a \in \sA$.

We now consider the equality  $\sB = {}\orth \smsub\orth \cap (\smsub[1])\orth \cap (\sA * \extn{\smsub})[1]$. 
We start by showing that $\sB \subseteq {}\orth \smsub\orth \cap (\smsub[1])\orth \cap (\sA * \extn{\smsub})[1]$. Let $b \in \sB = \rsms{\smsub}{\sA}$. Then there exists $a \in \sA$ and a minimal right $\extn{\smsub}$-approximation $\alpha \colon s_a \to a[1]$ such that $b$ occurs in a triangle, 
\[
s_a \rightlabel{\alpha} a[1] \too b \too s_a [1],
\]
so $\sB \subseteq  (\sA * \extn{\smsub})[1]$. By the triangulated Wakamatsu lemma (e.g. \cite[Lem.~2.1]{J09}), we have $b \in \smsub\orth$.

Applying $\Hom_{\sD}(-,\smsub)$ to the triangle above shows that $b \in {}\orth\smsub$. Indeed, the cases when $\smfull$ is an SMC or $w$-SMS with $w\geq 2$ follow from $\Hom_\sD (a[1], \smsub) = 0= \Hom_\sD (s_a[1],\smsub)$. The case when $\smfull$ is a $1$-SMS follows from Lemma~\ref{lem:Dugas-addendum}(2) using the fact that $\Hom_\sD (a[1], \smsub) = 0$. 

Applying $\Hom_{\sD}(\smsub,-)$ to this triangle gives an exact sequence,
\[
\begin{tikzcd}
\Hom_{\sD}(\smsub,a) \ar{r} & \Hom_{\sD}(\smsub,b[-1]) \ar{r} & \Hom_{\sD}(\smsub,s_a) \ar{rr}{\Hom_{\sD}(\smsub, \alpha)} & & \Hom_{\sD}(\smsub,a[1]),
\end{tikzcd}
\]
in which the morphism $\Hom_{\sD}(\smsub, \alpha)$ is an isomorphism by Lemma~\ref{lem:Dugas} and $\Hom_{\sD}(\smsub, a) = 0$ because $\sA \subseteq \smsub\orth$.
Hence, $\sB \subseteq {}\orth \smsub\orth \cap (\smsub[1])\orth \cap (\sA * \extn{\smsub})[1]$. 

Conversely, suppose $d \in {}\orth \smsub\orth \cap (\smsub[1])\orth \cap (\sA * \extn{\smsub})[1]$.
Since $d \in  (\sA * \extn{\smsub})[1]$, there is a triangle $s \rightlabel{\alpha} a[1] \too d \too s[1]$ with $a \in \sA$ and $s \in \extn{\smsub}$. We claim that $\alpha$ is a minimal right $\extn{\smsub}$-approximation from which it follows that $d \in \sB = \rsms{\smsub}{\sA}$, by definition.
As $d \in \smsub\orth$, it is immediate that $\alpha$ is a right $\extn{\smsub}$-approximation. Suppose $\alpha$ is not right minimal. It follows that $\alpha \colon s \to a[1]$ is isomorphic to 
$\begin{bmatrix} 
\alpha_1 & 0  
\end{bmatrix}
\colon s_1 \oplus s_2 \to a[1]$ 
with $s_1, s_2 \in \extn{\smsub}$.
In particular, $d$ has a direct summand isomorphic to $s_2 [1]$, contradicting the assumption that $d \in (\smsub[1])\orth$. Hence, $\alpha$ is a minimal right $\extn{\smsub}$-approximation.
It follows that $\sB \supseteq {}\orth \smsub\orth \cap (\smsub[1])\orth \cap (\sA * \extn{\smsub})[1]$. 

The proof that $(\lsms{\smsub}{\smfull \setminus \smsub}, \smfull \setminus \smsub)$ is an $\smsub$-mutation pair is similar.
\end{proof}

\begin{proposition} \label{prop:mutation-pair}
Let $\smsub$ be an orthogonal collection satisfying SMC Setup~\ref{SMC-setup} (SMS Setup~\ref{SMS-setup}, resp.).
Suppose $(\sA,\sB)$ is an $\smsub$-mutation pair. 
Then $\sA \cup \smsub$ is an SMC ($w$-SMS, resp.) in $\sD$ if and only if $\sB \cup \smsub$ is an SMC ($w$-SMS, resp.) in $\sD$.
\end{proposition}

\begin{proof}
Let $(\sA,\sB)$ be an $\smsub$-mutation pair. Suppose $\sA \cup \smsub$ is an SMC ($w$-SMS, resp.) in $\sD$. 
First observe that $\sA \subseteq \sZ$ because $\sA \cup \smsub$ is an SMC ($w$-SMS, resp.).
Therefore, $\sA$ is an SMC ($w$-SMS, resp.) in $\sZ$ by Theorem~\ref{thm:sm-reduction}.

If $\smsub$ satisfies SMC Setup~\ref{SMC-setup}, since $\sA \subseteq (\smsub[\geq 0])\orth$, then by Remark~\ref{rem:pair-implies-mutation} we can apply Lemma~\ref{lem:pair-implies-mutation}\eqref{right-approx-exists} to get $\sB = \rsms{\smsub}{\sA}$.
Similarly, applying Lemma~\ref{lem:pair-implies-mutation}\eqref{right-approx-exists} in the case that $\smsub$ satisfies SMS Setup~\ref{SMS-setup} gives $\sB = \rsms{\smsub}{\sA}$.
By Remark~\ref{rem:sm-reduction}, the right mutation triangle of Definition~\ref{def:mutation} 
coincides with the triangle used to define the shift in $\sZ$ in Theorem~\ref{thm:sm-reduction}.
It follows that $\sB = \sA\shift{1}$. In particular, $\sB \subseteq \sZ$, and since $\sA$ is an SMC ($w$-SMS, resp.) in $\sZ$, so is $\sB = \sA\shift{1}$. Applying Theorem~\ref{thm:sm-reduction} again, we obtain that $\sB \cup \smsub$ is an SMC ($w$-SMS, resp.) in $\sD$.

The proof of the other implication is similar, using Lemma~\ref{lem:pair-implies-mutation}\eqref{left-approx-exists}.
\end{proof}

This sufficiency of the SMC and SMS Setups in Theorem~\ref{thm:mut-SM-is-SM} now follows from Proposition~\ref{prop:mutation-pair}.
We now turn our attention to the necessity of the SMC and SMS Setups.

\subsection{Necessity of the SMC and SMS Setups} 

We start by recalling the following one-sided version of a result from \cite{Jin} relating SMC Setup~\ref{SMC-setup} with the existence of certain (unbounded) t-structures.

\begin{proposition}[{\cite[Prop. 3.2]{Jin}}] \label{prop:pre-smc-t-structure}
Let $\sD$ be a Hom-finite, Krull--Schmidt triangulated category and let $\smsub$ be an $\infty$-orthogonal collection in $\sD$. 
The pair $(\susp{}{\smsub}, (\smsub[\geq 0])\orth)$ is a t-structure in $\sD$ if and only if
\begin{enumerate}
\item the extension closure $\extn{\smsub}$ is contravariantly finite in $(\smsub[\geq 1])\orth$; and,
\item for any object $d$ in $\sD$ we have $\Hom_\sD(\smsub[\gg 0],d) = 0$.
\end{enumerate}
\end{proposition}

There is an obvious dual characterisation of when $({}\orth(\smsub[\leq 0], \cosusp{}{\smsub})$ is a t-structure.

\begin{proposition} \label{prop:smc-setup-is-necessary}
Let $\sD$ be a Hom-finite, Krull-Schmidt triangulated category and let $\smfull$ be an SMC in $\sD$. 
Suppose $\smsub \subseteq \smfull$. If $\rsmsn{n}{\smsub}{\smfull}$ is an SMC for each $n \geq 1$, then
\begin{enumerate}
\item the extension closure $\extn{\smsub}$ is contravariantly finite in $(\smsub[\geq 1])\orth$; and,
\item for any object $d$ in $\sD$ we have $\Hom_\sD(\smsub[\gg 0],d) = 0$.
\end{enumerate}
\end{proposition}

There is an obvious dual statement if $\lsmsn{n}{\smsub}{\smfull}$ is an SMC for each $n \geq 1$.

\begin{proof}
We first note that statement $(2)$ follows from $\smfull$ being an SMC; see Remark~\ref{rem:smc2-is-necessary}.

To establish statement $(1)$, we use Proposition~\ref{prop:pre-smc-t-structure} and Theorem~\ref{thm:length-simple-tilts}.
Write $\sH_0 = \extn{\smfull}$ for the heart of the bounded t-structure $(\susp{}{\smfull}, \cosusp{}{\smfull[-1]} = (\smfull[\geq 0])\orth)$.
We will show that $(\susp{}{\smsub[1]}, (\smsub[\geq 1])\orth)$ is a t-structure. Closure under shifts and the Hom-orthogonality are clear. It remains to find truncation triangles for each object $d$ of $\sD$.
Since $\smsub \subseteq \smfull$, we have $(\smfull[\geq 1])\orth \subseteq (\smsub[\geq 1])\orth$, and therefore it is sufficient to exhibit truncation triangles only for objects that lie in
\[
\susp{}{\smfull[1]} = \bigcup_{n \geq 1} \sH_0[n] * \cdots * \sH_0[2] * \sH_0[1].
\]

If $\rsmsn{i}{\smsub}{\smfull}$ is an SMC for each $i \geq 1$, then there is a sequence of length hearts and torsion pairs defined iteratively for $i \geq 0$ by simple right HRS tilts at the torsion classes $\cT_i = \extn{\smsub}[-i]$ in $\sH_i$:
\[
\sH_{i+1}  := \cF_i * \extn{\smsub}[-i-1], \text{ where }
\sH_i        = \extn{\smsub}[-i] * \cF_i \text{ with } \sH_0 := \extn{\smfull}.
\]
We proceed in three steps.

{\bf Step 1:} {\it For each $i \geq 0$, we have $\cF_i \subseteq (\smsub[\geq -i])\orth$.}

By definition of torsion pair, $\cF_i = \sH_i \cap \smsub[-i]\orth \subseteq \smsub[-i]\orth$. 
As $\sH_i$ is the heart of a (bounded) t-structure, we have $\Hom(\sH_i[\geq 1], \sH_i) = 0$, from which it follows that $\Hom(\smsub[\geq 1-i],\cF_i) = 0$ as $\smsub[-i] \subseteq \sH_i$. 

{\bf Step 2:} {\it For each $d \in \sH_0[n]$, $n \geq 1$, there exists a sequence of objects and morphisms
\[
d \to f_0[n] \to f_1[n] \to \cdots \to f_{n-1}[n]
\]
with $f_i \in \cF_i \subseteq \sH_i$ and such that $x := \cocone{d \to f_{n-1}[n]} \in \extn{\smsub}[n] * \cdots * \extn{\smsub}[2] * \extn{\smsub}[1]$.}

We proceed by induction on $n$. For $d \in \sH_0[1]$, the decomposition of $d[-1]$ with respect to the torsion pair $(\extn{\smsub},\cF_0)$ in $\sH_0$ gives the required triangle.

Now suppose $n > 1$ and $d \in \sH_0[n]$. By induction applied to $d[-1]$, there is a sequence of objects and morphisms
\[
d \to f_0[n] \to f_1[n] \to \cdots \to f_{n-2}[n]
\]
with $f_i \in \cF_i \subseteq \sH_i$ and such that $x_1 := \cocone{d \to f_{n-2}[n]} \in \extn{\smsub}[n] * \cdots * \extn{\smsub}[3] * \extn{\smsub}[2]$.
By construction, we have, 
\[
f_{n-2}[n] \in \cF_{n-2}[n] \subseteq \cF_{n-2}[n] * \extn{\smsub}[1] = \sH_{n-1}[n] = \extn{\smsub}[1] * \cF_{n-1}[n].
\]
Thus, we have a triangle $s_1[1] \to f_{n-2}[n] \to f_{n-1}[n] \to s_1[2]$. Applying the octahedral axiom to this triangle and the triangle $x_1 \to d \to f_{n-2}[n] \to x_1[1]$ gives the claim.

{\bf Step 3:} {\it For each $d \in \sH_0[n] * \cdots * \sH_0[2] * \sH_0[1]$ there exists a triangle
\[
x_d \to d \to y_d \to x_d[1],
\]
with $x_d \in \extn{\smsub}[n] * \cdots * \extn{\smsub}[2] * \extn{\smsub}[1]$ and $y_d \in \cF_{n-1}[n] * \cdots * \cF_1[2] * \cF_0[1]$.}

Again, we proceed by induction on $n$. For $n = 1$ the statement is contained in Step 2. 
Suppose $n > 1$ and $d \in \sH_0[n] * \cdots * \sH_0[2] * \sH_0[1]$. Decompose $d$ as $h[n] \to d \to c \to h[n+1]$ with $h \in \sH_0$ and $c \in \sH_0[n-1] * \cdots * \sH_0[2] * \sH_0[1]$.
By induction, there exists a triangle
\[
x_c \to c \to y_c \to x_c[1],
\]
with $x_c \in \extn{\smsub}[n-1] * \cdots * \extn{\smsub}[2] * \extn{\smsub}[1]$ and $y_c \in \cF_{n-2}[n-1] * \cdots * \cF_1[2] * \cF_0[1]$.
Applying the octahedral axiom twice, we get the following commutative diagrams:
\[
\begin{tikzcd}
                                  & x_c \ar[equals]{r} \ar{d} & x_c \ar{d} \\
d \ar[equals]{d} \ar{r} & c  \ar{r} \ar{d}                 & h[n+1] \ar{d} \\
d \ar{r}                       & y_c \ar{r}                        & x[1] 
\end{tikzcd}
\quad \text{and} \quad
\begin{tikzcd}
                                      & x_c[-1] \ar[equals]{r} \ar{d} & x_c[-1] \ar{d} \\
x_n \ar[equals]{d} \ar{r} & h[n]  \ar{r} \ar{d}                 & f_{n-1}[n] \ar{d} \\
x_n \ar{r}                       & x \ar{r}                        & y 
\end{tikzcd},
\]
where the middle horizontal triangle in the right-hand diagram comes from Step 2: in particular, $x_n \in \extn{\smsub}[n] * \cdots * \extn{\smsub}[1]$ and $f_{n-1} \in \cF_{n-1}$.
A final application of the octahedral axiom to the lower triangles in each of the two diagrams above gives us the following diagram.
\[
\begin{tikzcd}
                                  & y_c[-1] \ar[equals]{r} \ar{d} & y_c [-1] \ar{d} \\
x_n \ar[equals]{d} \ar{r} & x  \ar{r} \ar{d}                 & y \ar{d} \\
x_n \ar{r}                       & d \ar{r}                               & z
\end{tikzcd}
\]
We use $x_n \to d \to z \to x_n[1]$ to prove the existence of the required triangle.
Certainly, $x_n \in \extn{\smsub}[n] * \cdots * \extn{\smsub}[1]$. 
Reading from the right-most vertical triangles in the second and third octahedral axiom diagrams, we have
\begin{align*}
z & \in \big( \cF_{n-1}[n] * \extn{\smsub}[n-1] * \cdots *  \extn{\smsub}[1] \big) * \big( \cF_{n-2}[n-1] * \cdots * \cF_1[2] * \cF_0[1] \big) \\
       & \subseteq \big( \extn{\smsub}[n-1] * \cdots * \extn{\smsub}[3] * \cF_{n-1}[n] * \extn{\smsub}[2] * \extn{\smsub}[1] \big) * \big( \cF_{n-2}[n-1] * \cdots * \cF_0[1] \big) \\
       & \subseteq \big( \extn{\smsub}[n-1] * \cdots * \extn{\smsub}[2] * \cF_{n-1}[n] * \extn{\smsub}[1] \big) * \big( \cF_{n-2}[n-1] * \cdots * \cF_0[1] \big) \\
       & \subseteq \big( \extn{\smsub}[n-1] * \cdots * \extn{\smsub}[1] * \cF_{n-1}[n] \big) * \big( \cF_{n-2}[n-1] * \cdots * \cF_0[1] \big) \\
       & = \big( \extn{\smsub}[n-1] * \cdots * \extn{\smsub}[1] \big) * \big( \cF_{n-1}[n] * \cF_{n-2}[n-1] * \cdots * \cF_0[1] \big),
\end{align*}
where 
\begin{itemize}
\item the first inclusion follows from $\sH_{n-1} = \extn{\smsub}[1-n] * \cF_{n-1}$, so for $i = 3, \ldots, n-1$, we have that any triangle $\tilde{f}_{n-1}[n] \to e \to s[i] \to \tilde{f}_{n-1}[n+1]$  with $\tilde{f}_{n-1} \in \cF_{n-1}$ and $s \in \extn{\smsub}$ splits because $\Hom(\sH_{n-1}[\geq 1],\sH_{n-1}) = 0$ and $\extn{\smsub}[i] \subseteq \sH_{n-1}[n+i - 1]$;
\item the second inclusion follows from the fact that $(\extn{\smsub}[2], \cF_{n-1}[n+1])$ is a torsion pair in $\sH_{n-1}[n+1]$, so any triangle $\tilde{f}_{n-1}[n] \to e \to s[2] \to \tilde{f}_{n-1}[n+1]$  with $\tilde{f}_{n-1} \in \cF_{n-1}$ and $s \in \extn{\smsub}$ splits; and,
\item the third inclusion follows because $\cF_{n-1}[n] \subseteq \sH_{n-1}[n]$, $\extn{\smsub}[1] \subseteq \sH_{n-1}[n]$ and $\sH_{n-1}[n] = \extn{\smsub}[1] * \cF_{n-1}[n]$.
\end{itemize}
Hence, 
\begin{align*}
d & \in \big( \extn{\smsub}[n] * \cdots * \extn{\smsub}[1] \big) * \big( \extn{\smsub}[n-1] * \cdots * \extn{\smsub}[1] * \cF_{n-1}[n] * \cdots * \cF_0[1] \big) \\
   & \subseteq \big( \extn{\smsub}[n] * \cdots * \extn{\smsub}[1] \big) * \big( \cF_{n-1}[n] * \cdots * \cF_0[1] \big),
\end{align*}
where the inclusion follows again by \cite[Lem. 2.7]{CSP20} as $\Hom(\sH_0[\geq 1], \sH_0) = 0$. 

Thus, for any $d \in \sH_0[n] * \cdots *\sH_0[1]$, the claimed triangle, $x_d \to d \to y_d \to x_d[1]$, exists.
Moreover, by Step 1, we have $\cF_{n-1}[n] * \cdots * \cF_0[1] \subseteq (\smsub[\geq 1])\orth$, so this triangle is the truncation triangle for $d$ we were looking for.
Hence, $(\susp{}{\smsub[1]},(\smsub[\geq 1])\orth)$ is a t-structure and statements $(1)$ and $(2)$ hold by Proposition~\ref{prop:pre-smc-t-structure}. 
\end{proof}

Theorem~\ref{thm:mut-SM-is-SM} now follows; we briefly summarise below.

\begin{proof}[Proof of Theorem~\ref{thm:mut-SM-is-SM}]
($\Rightarrow$) for Statements (1) and (2): Apply Proposition~\ref{prop:mutation-pair} to the mutation pairs in Proposition~\ref{lem:mutation-pair-complement}, noting that the iteration of Proposition~\ref{lem:mutation-pair-complement} is allowed due to the fact that the conditions in the SMC Setup~\ref{SMC-setup}\eqref{smc1} are about $\sS$ and not about the SMC $\sS$ lies in.

($\Leftarrow$) for Statement (1): Applying Proposition~\ref{prop:smc-setup-is-necessary} and its dual gives conditions \eqref{smc1} and \eqref{smc2} of Setup~\ref{SMC-setup}.

($\Leftarrow$) for Statement (2): If $\sD$ is $(-w)$-Calabi--Yau, then any subcategory of $\sD$ is invariant under $\SSS[w]$. As $\smsub \subseteq \smfull$ and $\smfull$ is a $w$-SMS, $\extn{\smsub}$ is functorially finite in $\sD$ by \cite[Cor. 2.12]{CSP20}.
\end{proof}

\begin{remark}
Relaxing conditions SMC Setup~\ref{SMC-setup}\eqref{smc1} and \eqref{smc2} as in Remark~\ref{rem:sm-reduction}\ref{right-triangulated} and \ref{left-triangulated} permits infinite iteration of right $\smsub$-mutation or left $\smsub$-mutation because the reduced or quotient category $\sZ$ has a shift functor or loop functor, respectively. 
Applying Proposition~\ref{prop:smc-setup-is-necessary} and its dual shows that these conditions are again necessary and sufficient for right (resp. left) mutations of SMCs to be infinitely iterable.
\end{remark}

\subsection{Examples}

The following example illustrates the process in Theorem~\ref{thm:mut-SM-is-SM} for $w$-SMSs: reduce, shift, lift.

\begin{example}
Let $\sD = \Db (A_5)/\SSS [2]$, where $A_5$ is the linearly oriented Dynkin quiver of type $A_5$, and take $\smsub = \{s_1, s_2\}$ as in Figure~\ref{fig:typeA}. 

\begin{figure}
\centering
\begin{tikzpicture}[scale=0.8]

\draw[gray,rounded corners] (3,4.75) -- (0.25,-0.75) -- (10.75,-0.75) -- (8, 4.75) -- cycle;

\fill[red!10!white] (2,4) circle (2.5mm);
\draw[red,thick] (2,4) circle (2.5mm);

\fill[red!50!white,rounded corners] (1.7, 1.9) -- (2.1,1.7) -- (2.8, 3.1) -- (2.4, 3.3) -- cycle;
\draw[red,rounded corners,thick] (1.7, 1.9) -- (2.1,1.7) -- (2.8, 3.1) -- (2.4, 3.3) -- cycle;
\fill[red!50!white] (4,0) circle (2.5mm);
\draw[red,thick] (4,0) circle (2.5mm);

\fill[red!30!white,rounded corners] (4.7, 2.1) -- (5.1,2.3) -- (5.8, 0.9) -- (5.4, 0.7) -- cycle;
\draw[red,rounded corners,thick] (4.7, 2.1) -- (5.1,2.3) -- (5.8, 0.9) -- (5.4, 0.7) -- cycle;
\fill[red!30!white] (7,4) circle (2.5mm);
\draw[red,thick] (7,4) circle (2.5mm);

\fill[red!10!white,rounded corners] (7.7, 1.9) -- (8.1,1.7) -- (8.8, 3.1) -- (8.4, 3.3) -- cycle;
\draw[red,rounded corners,thick] (7.7, 1.9) -- (8.1,1.7) -- (8.8, 3.1) -- (8.4, 3.3) -- cycle;
\fill[red!10!white] (10,0) circle (2.5mm);
\draw[red,thick] (10,0) circle (2.5mm);

\fill[red!50!white,rounded corners] (9.7, 2.1) -- (10.1,2.3) -- (10.8, 0.9) -- (10.4, 0.7) -- cycle;
\draw[red,rounded corners,thick] (9.7, 2.1) -- (10.1,2.3) -- (10.8, 0.9) -- (10.4, 0.7) -- cycle;
\fill[red!50!white] (12,4) circle (2.5mm);
\draw[red,thick] (12,4) circle (2.5mm);

\fill[blue!30!white] (2,0) circle (2.5mm);
\draw[blue,thick] (2,0) circle (2.5mm);

\fill[blue!30!white,rounded corners] (4.5,2.6) -- (3.6,4.25) -- (5.4,4.25) -- cycle;
\draw[blue,rounded corners,thick] (4.5,2.6) -- (3.6,4.25) -- (5.4,4.25) -- cycle;

\fill[blue!30!white,rounded corners] (7.5,1.4) -- (6.6,-0.25) -- (8.4,-0.25) -- cycle;
\draw[blue,rounded corners,thick] (7.5,1.4) -- (6.6,-0.25) -- (8.4,-0.25) -- cycle;

\fill[blue!30!white] (10,4) circle (2.5mm);
\draw[blue,thick] (10,4) circle (2.5mm);

\fill[blue!30!white,rounded corners] (12.5,1.4) -- (11.6,-0.25) -- (13.4,-0.25) -- cycle;
\draw[blue,rounded corners,thick] (12.5,1.4) -- (11.6,-0.25) -- (13.4,-0.25) -- cycle;

\fill[blue!10!white] (5,0) circle (2.5mm);
\draw[blue,thick] (5,0) circle (2.5mm);

\fill[blue!10!white] (8,4) circle (2.5mm);
\draw[blue,thick] (8,4) circle (2.5mm);

\fill[blue!10!white] (13,4) circle (2.5mm);
\draw[blue,thick] (13,4) circle (2.5mm);

\node[above] (S) at (2.5,3.2) {\tiny $s$};
\node[below] (S0) at (10.5,0.8) {\tiny $s$};
\node[below] (S1) at (2,1.8) {\tiny $s_1$};
\node[below] (S2) at (4,-0.2) {\tiny $s_2$};
\node[above] (S1') at (10,2.2) {\tiny $s_1$};
\node[above] (S2') at (12,4.2) {\tiny $s_2$}; 
\node[below] (x1) at (7,-0.2) {\tiny $x_1$};
\node[below] (x2) at (8,-0.2) {\tiny $x_2$};
\node[above] (x3) at (8,4.2) {\tiny $x_3$};
\node[above] (y1) at (10,4.2) {\tiny $x_1\extn{1}$};
\node[below] (y2) at (2,-0.2) {\tiny $x_1\extn{1}$};
\node[above] (x[1]1) at (11,4.2) {\tiny $x_2[1]$};
\node[below] (x[1]2) at (3,-0.2) {\tiny $x_2[1]$};
\node[above] (x3[1]) at (3,4.2) {\tiny $x_3[1]$};
\node[below] (x3[1]1) at (11,-0.2) {\tiny $x_3[1]$};
\node[below] (x<1>1) at (4.5, 2.8) {\tiny $x_2\extn{1}$};
\node[above] (x<1>2) at (12.5, 1.2) {\tiny $x_2\extn{1}$};
\node[below] (t1) at (5, -0.2) {\tiny $x_3\extn{1}$};
\node[above] (t2) at (13,4.2) {\tiny $x_3\extn{1}$};

\foreach \x in {1,2,...,13}
	\foreach \y in {0,2,...,4}
{	
\fill (\x, \y) circle (0.5mm);
}
\foreach \x in {1.5,2.5,...,13}
	\foreach \y in {1,3}
{	
\fill (\x, \y) circle (0.5mm);
}
\end{tikzpicture}
\caption{Auslander-Reiten quiver of $\sD$ with arrows omitted. The part outlined in grey is a fundamental domain of $\sD$. The extension closure $\extn{\smsub}$ is shaded dark red, $\extn{\smsub}[1]$ mid-red, and $\extn{\smsub}[2]$ light red. $\sZ$ is shaded dark blue and light blue.}
\label{fig:typeA}
\end{figure}

Since $\sD$ has finitely many indecomposable objects (up to isomorphism), $\extn{\smsub}$ is functorially finite in $\sD$. Moreover, $\smsub$ is invariant under $\SSS[2]$ since $\sD$ is $(-2)$-CY. The subcategory $\sZ = \smsub^{\perp_{2}}$ of $\sD$ is the reduction of $\sD$ with respect to $\smsub$, and it is a $(-2)$-CY triangulated category. Moreover, we have
\[
\sZ \simeq \Db(A_2)/\SSS [2] \oplus \Db(A_1)/\SSS [2].
\]
The component $\Db(A_2)/\SSS [2]$ is indicated in dark blue in Figure~\ref{fig:typeA}, while the component $\Db(A_1)/\SSS [2]$ is indicated in light blue. Let $\sX = \{x_1, x_2, x_3\}$ as in Figure~\ref{fig:typeA}. The collection $\smfull = \smsub \cup \sX$ is a $2$-SMS in $\sD$, and $\sX = \smfull \setminus \smsub$ is a $2$-SMS in $\sZ$. Now, in order to do the right mutation of $\smfull$ with respect to $\smsub$, we compute the cones of the minimal right $\extn{\smsub}$-approximations $0 \to x_1[1], s_1 \to x_2[1]$ and $s \to x_3[3]$. Observe that the cones of these morphisms are $x_1\shift{1},  x_2\shift{1}$ and $x_3\shift{1}$, respectively, where $\shift{1}$ denotes the shift functor in $\sZ$. In other words, mutation of $\smfull$ with respect to $\smsub$ in $\sD$ corresponds to performing a shift of $\smfull \setminus \smsub$ in the reduction, $\sZ$.
\end{example}

The next example looks at SMC mutation outside the context of finite-dimensional algebras. We also use this example to illustrate the compatibility of SMC mutation with simple HRS tilting in Theorem~\ref{thm:length-simple-tilts}\ref{SMC} $\implies$ \ref{length}.
We refer also to \cite{Chen} for a more detailed study of SMC mutation in tube categories.

\begin{example} \label{ex:tube}
Let $\sH$ be a standard stable tube of rank $3$ and $\sD = \Db (\sH)$.
Consider $\smfull = \{s_1, s_2, s_3\}$, the set of simple objects on the mouth of the tube $\sH$, with $\tau s_i = s_{i-1}$, where $\tau$ denotes the Auslander--Reiten translate. Clearly, $\smfull$ is an SMC in $\sD$, with $\sH \coloneqq \extn{\smfull}$.

Now, let $\smsub = \{s_1, s_2\} \subseteq \smfull$. 
As $\cT \coloneqq \extn{\smsub} = \add{s_1, s_2, \begin{matrix} s_2\\[-2pt] s_1 \end{matrix}}$ has an additive generator, it is functorially finite in $\sD$. So in particular, $\extn{\smsub}$ satisfies condition (1) of SMC Setup~\ref{SMC-setup}. The collection $\smsub$ satisfies condition (2) by Remark~\ref{rem:smc2-is-necessary}.
Write $\cF = \smsub\orth \cap \sH$ for the corresponding torsionfree class of the torsion pair $\tt = (\cT,\cF)$.

In order to compute the right SMC mutation of $\smfull$ at $\smsub$, we compute the minimal right $\extn{\smsub}$-approximation of $s_3[1]$ and complete it to a triangle:
\[
\begin{matrix} s_2\\[-2pt] s_1 \end{matrix} \to s_3[1] \to \begin{matrix} s_2\\[-2pt] s_1\\[-2pt] s_3 \end{matrix}[1] \to \begin{matrix} s_2\\[-2pt] s_1 \end{matrix}[1].
\]
Thus $\rsms{\smsub}{\smfull} = \{s_1, s_2, \begin{matrix} s_2\\[-2pt] s_1\\[-2pt] s_3 \end{matrix}[1]\}$. Figure~\ref{fig:tube} illustrates the right simple tilt of $\sH$ at $\smsub$. We can see that $\rsms{\smsub}{\smfull}[-1]$ is indeed the set of simples of the new heart $\rhrs{\tt}{\sH}$. 
\begin{figure}
\begin{center}
\begin{tikzpicture}[scale=0.625]

\foreach \u in {0,5,10,15}{ 
\tikzmath{
\s = \u;
\q = \u + 0.5;
\p = \u + 2;
\r = \u + 2.5;
\t = \u + 3;}

\fill[red!80] (10,0) circle (1mm);
\fill[red!80] (11,0) circle (1mm);
\fill[red!80] (10.5,1) circle (1mm);
\fill[red!80] (13,0) circle (1mm);

\foreach \x in {0,1,2}
	{\fill[green!60!black] (12+0.5*\x,\x) circle (1mm);}
\foreach \x in {2,3,...,6}
	{\fill[green!60!black] (9+0.5*\x,\x) circle (1mm);}

\draw[dashed] (\s,0) -- (\s,6);
\draw[dashed] (\t,0) -- (\t,6);

 \foreach \y in {0,2,...,6} \foreach \x in {\s,...,\t}
  	{\fill[white] (\x,\y) circle (1mm);
   	\draw (\x,\y) circle (1mm);}
 \foreach \y in {1,3,5} \foreach \x in {\q,...,\r}
   	{\draw (\x,\y) circle (1mm);}
\foreach \y in {0,2,4}  \foreach \x in {\s,...,\p}
 	{ \node (node \y A) at (\x, \y) {};
 	\node (node \y B) at (\x+0.5, \y+1) {};
 	\draw[->] (node \y A) -- (node \y B);}
\foreach \y in {1,3,5}  \foreach \x in {\q,...,\r}
	 { \node (node \y A) at (\x, \y) {};
 	\node (node \y B) at (\x+0.5, \y+1) {};
 	\draw[->] (node \y A) -- (node \y B);}
\foreach \y in {2,4,6}  \foreach \x in {\s,...,\p}
 	{ \node (node \y A) at (\x, \y) {};
 	\node (node \y B) at (\x+0.5, \y-1) {};
 	\draw[->] (node \y A) -- (node \y B);}
\foreach \y in {1,3,5}  \foreach \x in {\q,...,\r}
	 { \node (node \y A) at (\x, \y) {};
 	\node (node \y B) at (\x+0.5, \y-1) {};
 	\draw[->] (node \y A) -- (node \y B);}
 }

\draw[red!80!black,thick] (10,0) circle (1mm);
\draw[red!80!black,thick] (11,0) circle (1mm);
\draw[green!30!black,thick] (12,0) circle (1mm);
\draw[red!80!black,thick] (13,0) circle (1mm);

\draw[decoration={brace,mirror},decorate] (10,-1) -- (13,-1);
\draw[decoration={brace},decorate] (10,7) -- (21,7);
\draw[decoration={brace},decorate] (-3,7) -- (8,7);

\node at (-2.5,3) {\small $\cdots$};
\node at (2.5,7.5) {\tiny $\sY$};
\node at (15.5,7.5) {\tiny $\sX$};
\node at (20.5,3) {\small $\cdots$};
\node at (10,-0.5) {\tiny $s_1$};
\node at (11,-0.5) {\tiny $s_2$};
\node at (12,-0.5) {\tiny $s_3$};
\node at (13,-0.5) {\tiny $s_1$};
\node at (11.5,-1.5) {\tiny $\sH$};

\end{tikzpicture}

\bigskip

\begin{tikzpicture}[scale=0.625]

\foreach \u in {0,5,10,15}{ 
\tikzmath{
\s = \u;
\q = \u + 0.5;
\p = \u + 2;
\r = \u + 2.5;
\t = \u + 3;}

\fill[red!50] (5,0) circle (1mm);
\fill[red!50] (6,0) circle (1mm);
\fill[red!50] (5.5,1) circle (1mm);
\fill[red!50] (8,0) circle (1mm);

\foreach \x in {0,1,2}
	{\fill[green!60!black] (12+0.5*\x,\x) circle (1mm);}
\foreach \x in {2,3,...,6}
	{\fill[green!60!black] (9+0.5*\x,\x) circle (1mm);}
	
\foreach \x in {0,1,2}
	{\fill[green!80!black] (7+0.5*\x,\x) circle (1mm);}
\foreach \x in {2,3,...,6}
	{\fill[green!80!black] (4+0.5*\x,\x) circle (1mm);}

\draw[dashed] (\s,0) -- (\s,6);
\draw[dashed] (\t,0) -- (\t,6);

 \foreach \y in {0,2,...,6} \foreach \x in {\s,...,\t}
  	{\fill[white] (\x,\y) circle (1mm);
   	\draw (\x,\y) circle (1mm);}
 \foreach \y in {1,3,5} \foreach \x in {\q,...,\r}
   	{\draw (\x,\y) circle (1mm);}
\foreach \y in {0,2,4}  \foreach \x in {\s,...,\p}
 	{ \node (node \y A) at (\x, \y) {};
 	\node (node \y B) at (\x+0.5, \y+1) {};
 	\draw[->] (node \y A) -- (node \y B);}
\foreach \y in {1,3,5}  \foreach \x in {\q,...,\r}
	 { \node (node \y A) at (\x, \y) {};
 	\node (node \y B) at (\x+0.5, \y+1) {};
 	\draw[->] (node \y A) -- (node \y B);}
\foreach \y in {2,4,6}  \foreach \x in {\s,...,\p}
 	{ \node (node \y A) at (\x, \y) {};
 	\node (node \y B) at (\x+0.5, \y-1) {};
 	\draw[->] (node \y A) -- (node \y B);}
\foreach \y in {1,3,5}  \foreach \x in {\q,...,\r}
	 { \node (node \y A) at (\x, \y) {};
 	\node (node \y B) at (\x+0.5, \y-1) {};
 	\draw[->] (node \y A) -- (node \y B);}
 }

\draw[red!80!black,thick] (5,0) circle (1mm);
\draw[red!80!black,thick] (6,0) circle (1mm);
\draw[green!30!black,thick] (10,2) circle (1mm);
\draw[green!30!black,thick] (13,2) circle (1mm);
\draw[red!80!black,thick] (8,0) circle (1mm);

\draw[decoration={brace,mirror},decorate] (5,-1) -- (13,-1);
\draw[decoration={brace},decorate] (10,7) -- (21,7);
\draw[decoration={brace},decorate] (-3,7) -- (3,7);

\node at (-2.5,3) {\small $\cdots$};
\node at (0,7.5) {\tiny $\sY[-1]$};
\node at (15.5,7.5) {\tiny $\sX$};
\node at (20.5,3) {\small $\cdots$};
\node at (5,-0.5) {\tiny $s_1'$};
\node at (6,-0.5) {\tiny $s_2'$};
\node at (8,-0.5) {\tiny $s_1'$};
\node at (9.5,2) {\tiny $s_3'$};
\node at (13.5,2) {\tiny $s_3'$};
\node at (9,-1.5) {\tiny $\rsms{\tt}{\sH} = \cF * \cT[-1]$};

\end{tikzpicture}
\end{center}
\caption{Top: the torsion pair $(\cT,\cF) = (\extn{s_1,s_2},(s_1,s_2)\orth \cap \sH)$ indicated inside $\sH =\extn{s_1, s_2, s_3}$ together with the bounded t-structure whose heart is $\sH$; 
$\cT$ is indicated in red, $\cF$ in green and the simple objects $\{s_1, s_2, s_3\}$ have thicker outlines.
Bottom: the negative shift of the torsion pair $(\cT,\cF)$ is shown in light red and green, respectively. The torsionfree class $\cF$ is shown in darker green and the right simple tilt $\rhrs{\tt}{\sH} = \cF * \cT[-1]$ is shown in light red and dark green with simple objects $\{s_1', s_2', s_3'\}$ shown with thicker outlines. 
The corresponding bounded t-structure is $(\sX',\sY') = (\sX * \cT[-1],\sY[-1] * \cF[-1])$.
Note that the non-zero morphisms and shifts go from left to right because $\sH$ is hereditary. 
\label{fig:tube}}
\end{figure}
\end{example}

\section{Silting and cosilting simple-minded collections} \label{sec:silting}

In this section, we will introduce the notions of silting, cosilting and bisilting SMCs in relation to the existence of adjacent co-t-structures.

Let $\Lambda$ be a finite-dimensional algebra. In \cite{Koenig-Yang}, Koenig and Yang show that right and left mutations of SMCs in $\Db(\Lambda)$ are always defined.  The key step in their proof is \cite[Lem.~7.8]{Koenig-Yang}, which establishes the existence of the $\extn{\smsub}$-approximations needed to define the mutation, cf. Theorem~\ref{thm:length-simple-tilts}\ref{T-approx}.
Their argument requires an involved construction using the realisation functor and the existence of enough projectives.

Our framework provides a conceptual homological understanding of why mutation of SMCs in $\Db(\Lambda)$ is always possible and how having enough projectives is sufficient for mutation of SMCs to be defined.
In addition, it also highlights that the difference between silting and cosilting objects can be identified in Hom-finite, Krull--Schmidt triangulated categories: for finite-dimensional algebras this is not detectable because silting objects are cosilting.

\subsection{Simple-minded collections and adjacent co-t-structures}

We define silting and cosilting SMCs via the existence of adjacent co-t-structures in the sense of \cite{Bondarko}.

\begin{definition} \label{def:silting-SMC}
Let $\smfull$ be an SMC and $(\sX, \sY) = (\susp{\sD}{\smfull}, \cosusp{\sD}{\smfull[-1]})$ the corresponding bounded t-structure. We say $\smfull$ is a:
\begin{enumerate}
\item \emph{silting SMC} if $(\sX, \sY)$ admits a left adjacent co-t-structure $({}\orth\sX, \sX)$.
\item \emph{cosilting SMC} if $(\sX, \sY)$ admits a right adjacent co-t-structure $(\sY, \sY\orth)$.
\item \emph{bisilting SMC} if it is both a silting and cosilting SMC. 
\end{enumerate}
The corresponding bounded t-structure $(\sX, \sY)$ will be called a \emph{silting t-structure}, \emph{cosilting t-structure} or \emph{bisilting t-structure}, respectively.
Similarly for its heart $\sH$, and we will use these terms interchangeably. 
\end{definition}

\begin{remark}\label{rem:strong-equivalent-conditions}
It follows from~\cite[Thm.~2.4]{CSPP} that if the t-structure $(\sX,\sY)$ admits a left adjacent co-t-structure then its heart $\sH$ satisfies the following:
\begin{enumerate}
\item $\sH$ is covariantly finite in $\sD$, and \label{covariantly-finite}
\item $\sH$ has enough projectives.  \label{enough-projs}
\end{enumerate}
There is a dual statement for cosilting SMCs.

In the case when $\sD$ is a saturated category (see \cite{Bondal-Kapranov}), e.g. $\sD=\Db(\Lambda)$ for a finite-dimensional algebra $\Lambda$ of finite global dimension or $\sD = \Db(\coh{X})$ for a smooth projective variety $X$, then \eqref{covariantly-finite} is equivalent to \eqref{enough-projs}, which is equivalent to the existence of a left adjacent co-t-structure by \cite[Cor.~2.8]{CSPP}.
In particular, for a saturated triangulated category we have
\begin{align*}
\text{$(\sX,\sY)$ is bisilting} & \iff  \text{$\sH$ is functorially finite in $\sD$} \\
 & \iff \text{$\sH$ has enough projectives and enough injectives.}
 \end{align*}
\end{remark}

\begin{remark} \label{rem:silting-SMC}
Remark~\ref{rem:strong-equivalent-conditions} explains our choice of terminology in Definition~\ref{def:silting-SMC}. 
Indeed, in the case that $(\sX,\sY)$ admits a left adjacent t-structure, 
in the language of \cite[\S 2]{CSPP}, \cite[Thm. 2.4]{CSPP} says that each projective object $q$ of $\sH$ admits an $\sH$-epimorphism $H^0(p) \onto q$ for some object $p$ of the \emph{projective coheart} $\sP = {}\orth \sX[1] \cap \sX$, which is a presilting subcategory, i.e.\ $\Hom(\sP,\sP[>0])=0$ (see \cite{AI}).
Dually, when $(\sX, \sY)$ admits a right adjacent co-t-structure, each injective object $e$ of $\sH$ admits an $\sH$-monomorphism $e \into H^0(i)$ for some object $i$ of the \emph{injective coheart} $\sI = \sY\orth \cap \sY[1]$, which again is a presilting subcategory.

In the first case, the coheart of the left adjacent co-t-structure $\sP$ is `projective minded' because it `approximates' the projective objects of $\sH$, so we call the SMC generating $\sH$ `silting'.
In the second case, the coheart of the right adjacent co-t-structure $\sI$ is `injective minded' because it `approximates' the injective objects of $\sH$, so we call the SMC generating $\sH$ `cosilting'.
In the `bisilting' case, there are two adjacent co-t-structures, one whose coheart is projective minded and one whose coheart is injective minded.
\end{remark}

\begin{remark}\label{rem:adjacent-co-t-structures}
Let $(\sX, \sY)$ be a t-structure. In particular, $\sX$ is contravariantly finite in $\sD$ and $\sY$ is covariantly finite in $\sD$. The t-structure $(\sX, \sY)$ admits a left (resp. right) adjacent co-t-structure if and only if $\sX$ is covariantly finite in $\sD$ (resp. $\sY$ is contravariantly finite in $\sD$). Therefore $(\sX,\sY)$ is a bisilting t-structure if and only if $\sX$ and $\sY$ are both functorially finite in $\sD$. 
\end{remark}

\begin{example}
Recall the tube category $\sD = \Db (\sH)$ from Example~\ref{ex:tube}, where $\sH$ is a standard, stable tube of rank $3$. Recall that $\smfull=\{s_1,s_2,s_3\}$ is an SMC with $\extn{\smfull} = \sH$. However, $\extn{\smfull}$ is neither covariantly finite nor contravariantly finite in $\sD$, and so by Remark~\ref{rem:strong-equivalent-conditions}, $\smfull$ is neither a silting nor a cosilting SMC.
\end{example}

\begin{example} \label{ex:silting-not-cosilting}
Consider the following quiver of type $A_\infty$:
\[
Q\colon \cdots \to 3 \to 2 \to 1,
\] 
and let $\sD = \Db (\repb{Q})$, where $\repb{Q}$ denotes the category of finite dimensional representations of $Q$. We have that $\repb{Q}$ is an hereditary abelian category with enough projectives but not enough injectives (see~\cite[Props.~1.15 \& 1.16]{BLP}). Moreover, the projective objects in $\repb{Q}$ are of the form $P_x = \begin{matrix} S_x\\[-2pt] \vdots \\[-2pt] S_1 \end{matrix}$, with $x \geq 1$. Denote by $\cP$ the set of projective objects in $\repb{Q}$. 

Take $\smfull = \{S_i \mid i \in Q_0\} \subseteq \repb{Q}$. Since $\smfull$ is $\infty$-orthogonal and $\extn{\smfull} = \repb{Q}$, it follows that $\smfull$ is an SMC in $\sD$. But $\extn{\smfull}$ doesn't have enough injectives, so $\smfull$ is not a cosilting SMC by Remark~\ref{rem:strong-equivalent-conditions}. On the other hand, the two conditions on the projectives of $\extn{\smfull}$ in~\cite[Thm.~2.4 (2)]{CSPP} are satisfied. Indeed, as already mentioned $\extn{\smfull}$ has enough projectives and one can check that the projective coheart $\,^\perp \sX[1] \cap \sX$, where $\sX = \susp{\sD}{\smfull}$, is $\cP$. It thus follows by~\cite[Thm.~2.4]{CSPP} that $\smfull$ is a silting SMC. 

Taking the same collection of objects $\smfull$ over the opposite quiver $Q^{op}$ would give a cosilting SMC which is not silting. 
\end{example}

\begin{proposition} \label{prop:hrs-preserves-strongly-algebraic}
Let $(\sX, \sY)$ be a t-structure with heart $\sH$ and let $\tt = (\cT,\cF)$ be a torsion pair in $\sH$.
\begin{enumerate}
\item Suppose $(\sX,\sY)$ admits a right adjacent co-t-structure and suppose that the torsionfree class $\cF$ is contravariantly finite in $\sH$.   Then the right HRS-tilt of $(\sX, \sY)$ at $\tt$ also admits a right adjacent co-t-structure.
\item Suppose $(\sX,\sY)$ admits a left adjacent co-t-structure and suppose that the torsion class $\cT$ is covariantly finite in $\sH$.   Then the left HRS-tilt of $(\sX, \sY)$ at $\tt$ also admits a left adjacent co-t-structure.
\end{enumerate}
\end{proposition}

\begin{proof}
We prove statement (1); statement (2) is dual.

The right HRS-tilt is $(\sX * \cT[-1], (\cF * \sY)[-1])$. As observed in Remark~\ref{rem:adjacent-co-t-structures}, it is enough to show that $\cF * \sY$ is contravariantly finite in $\sD$. By hypothesis, $\cF$ is contravariantly finite in $\sH$ and, by Remark~\ref{rem:strong-equivalent-conditions}, $\sH$ is contravariantly finite in $\sD$. Therefore, $\cF$ is contravariantly finite in $\sD$. Thus, given $c \in \sD$, we can then take a right $\cF$-approximation $\alpha \colon f \to c$ and extend it to a triangle $f \rightlabel{\alpha} c \too b \too f[1]$. 

Since $\sY$ is contravariantly finite in $\sD$, there is a right $\sY$-approximation $\beta \colon y \to b$ which extends to a triangle $y \rightlabel{\beta} b \too d \too y[1]$.
Applying the octahedral axiom we get the following commutative diagram.
\[
\begin{tikzcd}
f \ar{r} \ar[equals]{d} & y' \ar{r} \ar{d}{\gamma} & y \ar{d}{\beta} \\
f \ar{r}[swap]{\alpha} & c \ar{r} \ar{d}                 & b \ar{d} \\
                                  & d \ar[equals]{r}             & d
\end{tikzcd}
\]
It now follows by the dual of \cite[Lem.~5.3]{Saorin-Zvonareva} that $\gamma \colon y' \to c$ is a right $(\cF * \sY)$-approximation.
Hence, $\cF * \sY$ is contravariantly finite in $\sD$.
\end{proof}

\subsection{Finite bisilting SMCs are preserved under mutation}

It is natural to ask if mutation of a bisilting SMC is always defined and whether the bisilting property is preserved.
The main result of this section asserts that this is the case for finite SMCs.

\begin{theorem} \label{thm:bisilting-smc-mutation}
Let $\sD$ be a Hom-finite, Krull--Schmidt, $\kk$-linear triangulated category. 
Suppose $\smfull$ is a finite SMC and $\smsub \subseteq \smfull$.
\begin{enumerate}
\item If $\smfull$ is cosilting then $\rsmsn{n}{\smsub}{\smfull}$ is also a cosilting SMC for each $n \in \bN$.
\item If $\smfull$ is silting then $\lsmsn{n}{\smsub}{\smfull}$ is also a silting SMC for each $n \in \bN$.
\end{enumerate}
\end{theorem}

We note that Theorem~\ref{thm:bisilting-smc-mutation} can be reformulated in the language of simple tilts:

\begin{theorem} \label{thm:strong-hrs-tilt}
Let $\sH$ be the heart of a t-structure in $\sD$ which is length with finitely many simple objects. Suppose $\smsub$ is a subset of the simple objects of $\sH$. 
\begin{enumerate}
\item If $\sH$ is cosilting, then the right HRS-tilt of $\sH$ at the torsion pair $\tt = (\extn{\smsub},\smsub\orth\cap \sH)$ is again cosilting.
\item If $\sH$ is silting, then the left HRS-tilt of $\sH$ at the torsion pair $\tt = ({}\orth\smsub \cap \sH, \extn{\smsub})$ is again silting.
\end{enumerate}
\end{theorem}

We start by observing that appropriate one-sided mutations of silting or cosilting SMCs are SMCs.

\begin{lemma} \label{lem:smc-mutation}
Let $\sD$ be a Hom-finite, Krull--Schmidt, $\kk$-linear triangulated category. Let $\smfull$ be a (not necessarily finite) SMC and $\smsub \subseteq \smfull$. 
\begin{enumerate}
\item If $\smfull$ is cosilting, then $\rsmsn{n}{\smsub}{\smfull}$ is an SMC for all $n \in \bN$.
\item If $\smfull$ is silting, then $\lsmsn{n}{\smsub}{\smfull}$ is an SMC for all $n \in \bN$.
\end{enumerate}
\end{lemma}

\begin{proof}
We show (1); (2) is dual. 
To see that $\rsmsn{n}{\smsub}{\smfull}$ is an SMC, it is sufficient to check that the relevant half of SMC Setup~\ref{SMC-setup} holds. We can then apply Theorem~\ref{thm:mut-SM-is-SM}.

By Remark~\ref{rem:strong-equivalent-conditions}, we have that $\extn{\smfull}$ is contravariantly finite in $\sD$. Furthermore, by Lemma~\ref{lem:simples-are-funct-finite}, $\extn{\smsub}$ is functorially finite in $\extn{\smfull}$. It thus follows that $\extn{\smsub}$ is contravariantly finite in $\sD$. In particular, (the relevant half of) condition  \eqref{smc1} of Setup~\ref{SMC-setup} is satisfied. Condition \eqref{smc2} of Setup~\ref{SMC-setup} holds since $\smsub \subseteq \smfull$ and the condition holds for $\smfull$ by Remark~\ref{rem:smc2-is-necessary}.
\end{proof}

We now proceed to show that the bisilting property is preserved. For this we use simple HRS tilting and its compatibility with mutation; see the implication \ref{SMC} $\implies$ \ref{length} in Theorem~\ref{thm:length-simple-tilts}. We need the following lemmas.

\begin{lemma}[{\cite[Thm.]{Smalo}}] \label{lem:funct-finite}
Let $\Lambda$ be an artin algebra and $(\cT,\cF)$ a torsion pair in $\mod{\Lambda}$. The torsion class $\cT$ is functorially finite in $\mod{\Lambda}$ if and only if the torsionfree class $\cF$ is functorially finite in $\mod{\Lambda}$.
\end{lemma}

\begin{lemma}[{\cite[Lem.~6]{Al-Nofayee}}] \label{lem:module-category}
Let $\sH$ be a $\kk$-linear, Hom-finite abelian category. The following statements are equivalent.
\begin{enumerate}[label=(\roman*)]
\item There is an equivalence of categories $\sH \simeq \mod{\Lambda}$, where $\Lambda$ is a finite-dimensional $\kk$-algebra.
\item The category $\sH$ has a projective generator.
\item The category $\sH$ has an injective cogenerator.
\end{enumerate}
\end{lemma}

\begin{corollary}[cf. {\cite[Cor.~2.11]{CSPP}}] \label{cor:module-category}
Let $\sD$ be a Hom-finite, Krull--Schmidt, $\kk$-linear triangulated category. Suppose $\smfull$ is a finite SMC in $\sD$. If $\smfull$ is also silting or cosilting then $\extn{\smfull} \simeq \mod{\Lambda}$, where $\Lambda$ is a finite-dimensional $\kk$-algebra.
\end{corollary}

\begin{proof}
We show that if $\smfull$ is a finite cosilting SMC then $\extn{\smfull} \simeq \mod{\Lambda}$ for a finite-dimensional $\kk$-algebra $\Lambda$; the case when $\smfull$ is a finite silting SMC is analogous.

As $\smfull$ is cosilting, $\extn{\smfull}$ has enough injectives by Remark~\ref{rem:strong-equivalent-conditions}. 
It follows that the direct sum of the injective envelopes of each of the finitely many simple objects (since $\smfull$ is finite) gives an injective cogenerator of $\extn{\smfull}$. 
Since $\extn{\smfull}$ is a $\kk$-linear, Hom-finite abelian category, Lemma~\ref{lem:module-category} gives $\extn{\smfull} \simeq \mod{\Lambda}$ for some finite-dimensional $\kk$-algebra $\Lambda$.
\end{proof}

We are now ready to prove the main result of this section.

\begin{proof}[Proof of Theorem~\ref{thm:bisilting-smc-mutation}]
We show the first statement; the second is dual.
Suppose $\smfull$ is a finite cosilting SMC in $\sD$.
By Lemma~\ref{lem:smc-mutation}, $\rsms{\smsub}{\smfull}[-1]$ is a finite SMC. We therefore only need to check that $\rsms{\smsub}{\smfull}[-1]$ is cosilting.
By definition, $\sH = \extn{\smfull}$ is the heart of a cosilting t-structure $(\sX,\sY)$ and $\cT \coloneqq \extn{\smsub}$ is a functorially finite torsion class with torsionfree class $\cF = \smsub\orth \cap \sH$. By Corollary~\ref{cor:module-category}, $\sH \simeq \mod{\Lambda}$ for a finite-dimensional $\kk$-algebra $\Lambda$. Hence, by Lemma~\ref{lem:funct-finite}, the torsionfree class $\cF$ is also functorially finite.
By Proposition~\ref{prop:hrs-preserves-strongly-algebraic}, the right HRS-tilt, $\rhrs{\tt}{\sH}$ at $\tt = (\cT,\cF)$ is cosilting. Moreover, by the proof of~\ref{SMC} $\implies$ \ref{length} in Theorem~\ref{thm:length-simple-tilts}, we have $\rhrs{\tt}{\sH} = \extn{\rsms{\smsub}{\smfull}[-1]}$. That is, $\rsms{\smsub}{\smfull}[-1]$ is a cosilting SMC.
\end{proof}

\subsection{The case of module categories} \label{sec:mod-cat}

In this section we show how the Koenig--Yang correspondences in \cite{Koenig-Yang} can be interpreted as a stronger version of Theorem~\ref{thm:bisilting-smc-mutation} in the case that $\sD = \Db(\Lambda)$ for a finite-dimensional $\kk$-algebra $\Lambda$.
Here we see that the notions of silting SMCs and cosilting SMCs coincide, i.e. that all SMCs are bisilting.
This framework helps give a conceptual homological explanation of the compatibility of silting and SMC mutation observed in \cite{Koenig-Yang}. 

\begin{proposition} \label{prop:KY}
Let $\Lambda$ be a finite-dimensional $\kk$-algebra and $\sD = \Db (\Lambda)$. Every SMC in $\sD$ is a finite bisilting SMC. 
\end{proposition}

The key to this observation is the following lemma about the existence of co-t-structures, cf. Proposition~\ref{prop:pre-smc-t-structure}.
Recall that a subcategory $\sP$ of $\sD$ is \emph{presilting} if $\Hom_\sD (\sP, \sP[>0]) = 0$.

\begin{lemma}[{\cite[Prop.~3.2]{IY}}] \label{lem:IY-adj-co-t-structure}
Let $\sD$ be a triangulated category and $\sP$ be a presilting subcategory of $\sD$.
\begin{enumerate}[label=(\arabic*)]
\item \label{left-adjacent} The pair $(\cosusp{\sD}{\sP[-1]}, (\sP[<0])\orth)$ is a co-t-structure in $\sD$ if and only if the following conditions hold:
\begin{enumerate}
\item[(P1)] $\sP$ is contravariantly finite in $(P[<0])\orth$; and,
\item[(P2)] $\Hom (\sP[\ll 0], d) = 0$, for all $d \in \sD$. 
\end{enumerate}
\item \label{right-adjacent} The pair $({}\orth(\sP[\geq 0]), \susp{\sD}{\sP})$ is a co-t-structure in $\sD$ if and only if the following conditions hold:
\begin{enumerate}
\item[(P1$'$)]  $\sP$ is covariantly finite in ${}\orth(\sP[> 0])$; and,
\item[(P2$'$)] $\Hom (d,\sP[\gg 0]) = 0$, for all $d \in \sD$. 
\end{enumerate}
\end{enumerate}
In each case the coheart is $\sP$.
\end{lemma}

\begin{proof}[Proof of Proposition~\ref{prop:KY}]
It is clear that any SMC in $\sD = \Db (\Lambda)$ is finite, with cardinality  the rank of $K_0 (\sD)$.
Suppose $\smfull$ is an SMC in $\sD$ and let $(\sX, \sY) = (\susp{\sD}{\smfull}, \cosusp{\sD}{\smfull[-1]})$ be the corresponding bounded t-structure in $\sD$. To see that $\smfull$ is silting, we need to show that $(\sX, \sY)$ has a left adjacent co-t-structure. 
By \cite[Thm.~6.1]{Koenig-Yang}, 
\[
(\sX, \sY) = ((\sP[<0])\orth, (\sP[\geq 0])\orth),
\]
 where $\sP = {}^\perp \sX[1] \cap \sX$ is a silting subcategory of $\Kb(\proj{\Lambda})$. 
By \cite[Prop.~2.20]{AI}, $\sP$ has an additive generator. Hence, $\sP$ is a presilting subcategory of $\sD$ satisfying condition (P1) above. 
Now, since $\sP \subseteq \Kb(\proj{\Lambda})$, we have $\Hom (\sP[\ll 0],d)= 0$ for all $d \in \sD$ because $\sD \simeq \sK^{b,-} (\proj{\Lambda})$.  Thus, (P2) is also satisfied and, by Lemma~\ref{lem:IY-adj-co-t-structure}\ref{left-adjacent}, $(\sX,\sY)$ admits a left adjacent co-t-structure and $\smfull$ is silting.

Dually, the subcategory $\sI = \sY^\perp \cap  \sY[1]$ is silting in $\Kb (\inj{\Lambda})$. By the dual of \cite[Thm.~6.1]{Koenig-Yang} we have $(\sX, \sY) = ({}\orth(\sI[<0]), {}\orth(\sI[\geq 0]))$. It follows from Lemma~\ref{lem:IY-adj-co-t-structure}\ref{right-adjacent}, that $\smfull$ is a cosilting SMC if and only if conditions (P1$'$) and (P2$'$) are satisfied. Condition (P1$'$) follows from the fact that $\sI$ has an additive generator, and condition (P2$'$) from the fact that $\sI \subseteq K^b (\inj{\Lambda})$ and $\sD \simeq \sK^{b,+} (\inj{\Lambda})$. 
\end{proof}

\begin{remark}
We elaborate on the discussion in Remark~\ref{rem:silting-SMC}.
In the context of `small' (i.e.\ Hom-finite, Krull--Schmidt) triangulated categories, the notion `cosilting subcategory' has not really appeared.
The distinction between silting subcategories and cosilting subcategories typically occurs for `large' (i.e. admitting set-indexed products and coproducts) triangulated categories, where their definitions are distinguished by their behaviour with respect to coproducts and products, respectively; see \cite{ALSV} for an overview.

Philosophically, the distinction is determined by the side on which the adjacent t-structure occurs, meaning that even in the context of small triangulated categories there is a distinction between silting and cosilting subcategories, detected by the corresponding SMC via adjacency, despite the fact that finite coproducts and finite products coincide.

For finite-dimensional algebras of finite global dimension, the equivalences $\sK^b(\proj{\Lambda}) \simeq \Db(\Lambda) \simeq \sK^b(\inj{\Lambda})$ mean that one can view silting objects as either projective-minded objects or injective-minded objects. 
\end{remark}

Proposition~\ref{prop:KY} can be extended to the context of saturated triangulated categories. Here, the notions of silting and cosilting coincide for finite SMCs.

\begin{proposition}
Let $\sD$ be a Hom-finite, Krull--Schmidt, saturated triangulated category. Then a finite SMC is silting if and only if it is cosilting.
\end{proposition}

\begin{proof}
Let $\smfull$ be a finite SMC in $\sD$ with associated bounded t-structure $(\sX,\sY)$ and heart $\sH = \extn{\smfull}$. Then
\begin{align*}
\smfull \text{ is silting } & \iff (\sX,\sY) \text{ has a left adjacent co-t-structure} \\
                               & \iff \sH \text{ has enough projectives by \cite[Cor. 2.8]{CSPP}} \\
                               & \iff \sH \simeq \mod{\Lambda} \text{ for some finite-dimensional algebra } \Lambda \text{ by Lemma~\ref{lem:module-category}} \\
                               & \iff \sH \text{ has enough injectives by Lemma~\ref{lem:module-category}} \\
                               & \iff \sH \text{ has a right adjacent co-t-structure by \cite[Cor. 2.8]{CSPP}} \\
                               & \iff \smfull \text{ is cosilting,}
\end{align*}
where the third and fourth bi-implications use that $\smfull$ is finite to obtain a projective generator from the projective covers of the finitely many simple objects, and an injective cogenerator from the injective envelopes of the finitely many simple objects, respectively.
\end{proof}

Note that finiteness of $\smfull$ is crucial in the argument above: Example~\ref{ex:silting-not-cosilting} gives an example of an infinite SMC that is silting but not cosilting, which can be dualised to give an example of an infinite SMC that is cosilting but not silting.
We end with the following question.

\begin{question}
Let $\sD$ be a Hom-finite, Krull--Schmidt triangulated category. Is it always true that a finite SMC is silting if and only if it is cosilting?
\end{question}

We note that we expect the answer to this question to be false if $\sD$ is not Hom-finite and Krull--Schmidt.

\section{SMC mutations vs SMS mutations}

In \cite{Jin}, Haibo Jin established a relationship between SMC reduction and SMS reduction via a singularity category construction. In this section we observe that mutation of finite bisilting SMCs is compatible with mutation of SMSs via this construction. 

Before stating the result, we need to recall the definition of CY-triple  \cite[Def.~4.1]{Jin}. 

\begin{definition}
Let $w \geq 1$. A \emph{$(1-w)$-CY triple} is a tuple $(\sD, \sD^p, \smfull)$ where: 
\begin{enumerate}
\item $\sD$ is a Hom-finite, Krull--Schmidt, $\kk$-linear triangulated category and $\sD^p$ is a thick subcategory of $\sD$;
\item The functor $[1-w]$ is a {\it relative Serre functor}, i.e.\ it satisfies the bifunctorial isomorphism $\Hom (x,y) \cong D\Hom (y, x[1-w])$ for any $x \in \sD^p$ and $y \in \sD$; and,
\item $\smfull$ is a bisilting SMC in $\sD$ for which ${}\orth(\smfull[\geq 0])$ and $(\smfull[<0])\orth$ are subcategories of $\sD^p$. 
\end{enumerate}
\end{definition}

Given a $(1-w)$-CY triple, the {\it singularity category} $\Dsg$ is defined to be the Verdier quotient $\sD/\sD^p$. This category is a $(-w)$-CY triangulated category \cite[Thm.~4.5]{Jin}. 
The canonical quotient functor is denoted $\pi\colon \sD \to \Dsg$.
The prototypical example is the following.

\begin{example}
Let $\Lambda$ be a finite-dimensional symmetric algebra and write $\smfull$ for the set of simple $\Lambda$-modules. Then $(\Db(\Lambda), \Kb(\proj{\Lambda}), \smfull)$ is a $0$-CY triple. The singularity category $\Dsg = \Db(\Lambda)/\Kb(\proj{\Lambda})$ is the classic singularity category, which is equivalent to $\stmod{\Lambda}$ by a famous result of Ragnar-Olaf Buchweitz \cite{Buchweitz}.
\end{example}

\begin{theorem} \label{thm:compatibility}
Let $(\sD, \sD^p, \smfull)$ be a $(1-w)$-CY triple in which the bisilting SMC, $\smfull$, is finite. Let $\smsub$ be an $\infty$-orthogonal collection in $\sD$ such that $\extn{\smsub}$ is functorially finite in $\sD$. The following diagram is commutative:
\[
\begin{tikzcd}
\{\text{bisilting SMCs in $\sD$ containing $\smsub$}\} \ar{r}{\rsms{\smsub}{-}} \ar{d}{\pi} & \{\text{bisilting SMCs in $\sD$ containing $\smsub$}\} \ar{d}{\pi}\\
\{\text{$w$-SMSs in $\Dsg$ containing $\pi(\smsub)$}\} \ar{r}{\rsms{\pi(\smsub)}{-}}  & \{\text{$w$-SMSs in $\Dsg$ containing $\pi(\smsub)$}\}.
\end{tikzcd}
\]
\end{theorem}
\begin{proof}
We first note that the maps are well defined. Indeed, since $\sD$ has a finite SMC by assumption, every SMC in $\sD$ is finite. The top horizontal map is well defined by Theorem~\ref{thm:bisilting-smc-mutation}. 
The bottom horizontal map is well defined because $\sD_{sg}$ is $(-w)$-CY and functorial finiteness of $\extn{\pi (\smsub)}$ in $\sD_{sg}$ follows automatically if $\pi(\smsub)$ is a subset of a $w$-SMS by \cite[Cor.~2.12]{CSP20}.

To see that the vertical map is well defined, we must verify the hypothesis of \cite[Thms.~4.5 \& 4.13]{Jin}, the application of which then gives the claim. That is, we must check, given a finite bisilting SMC, $\sV$, in $\sD$ containing $\smsub$, 
\begin{enumerate}[label=(\roman*)]
\item its extension closure $\extn{\sV}$ is functorially finite in $\sD$, and, \label{funct-finite}
\item there exists $n \geq 0$ such that $\sV \subseteq \extn{\smfull}[n] * \extn{\smfull}[n-1] * \cdots * \extn{\smfull}[1-n] * \extn{\smfull}[-n]$. \label{bounded} 
\end{enumerate}
As $\sV$ is finite and $\extn{\smfull}$ is the heart of a bounded t-structure, hypothesis \ref{bounded} holds. Hypothesis \ref{funct-finite} follows from the assumption that $\sV$ is bisilting by Remark~\ref{rem:strong-equivalent-conditions}. 

Consider the diagram below, the back face of which is our desired commutative diagram. It suffices to show the remaining faces commute.
\begin{center}
\begin{tikzcd}[column sep=tiny]

\left\{ \text{\parbox{3.25cm}{\centering Bisilting SMCs in \\ $\sD$ containing $\sS$}} \right\} 
\ar{dd}[near start]{\pi}
\ar[<->]{rr}[near end]{\rsms{\sS}{-}}[near end,swap]{\sim} 
\ar[<->]{dr}{\alpha}[swap]{\sim} 
&& 
\left\{ \text{\parbox{3.25cm}{\centering Bisilting SMCs in \\ $\sD$ containing $\sS$}} \right\} 
\ar{dd}[near start]{\pi} 
\ar[<->]{dr}{\alpha}[swap]{\sim}
& \\
& \left\{ \text{\parbox{2cm}{\centering Bisilting SMCs in $\sZ$}} \right\} 
\ar[<->, crossing over]{rr}[near start]{\shift{1}_\sZ}[swap, near start]{\sim}
&& 
\left\{ \text{\parbox{2cm}{\centering Bisilting SMCs in $\sZ$}} \right\} 
\ar{dd}[near start]{\pi_\sZ}
\\
\left\{ \text{\parbox{3cm}{\centering $w$-SMSs in $\Dsg$ containing $\pi(s)$}} \right\} 
\ar[<->]{dr}{\alpha_{sg}}[swap]{\sim}
\ar[<->]{rr}[near end]{\rsms{\sS}{-}}[swap, near end]{\sim} 
 &&
\left\{ \text{\parbox{3.25cm}{\centering $w$-SMSs in $\Dsg$ containing $\pi(\sS)$}} \right\} 
\ar[<->]{dr}{\alpha_{sg}}[swap]{\sim}
\\
& \left\{ \text{\parbox{2cm}{\centering $w$-SMSs in $\Zsg$}} \right\} 
\ar[from=uu, crossing over]{}[near start]{\pi_\sZ} 
\ar[<->]{rr}[near start]{ \shift{1}_{sg}}[near start, swap]{\sim}
&&
\left\{ \text{\parbox{2cm}{\centering $w$-SMSs in $\Zsg$}} \right\} 
\end{tikzcd}
\end{center}

The top and bottom faces commute by Proposition~\ref{prop:mutation-pair}. 
The bijection between SMCs in $\sD$ containing $\smsub$ and SMCs in the reduction $\sZ$ of $\sD$ with respect to $\smsub$ reduces to a bijection between bisilting SMCs in $\sD$ containing $\smsub$ and bisilting SMCs in $\sZ$ by \cite[Thm.~6.1]{Jin}. The commutativity of the left and right faces follows from \cite[p.~1485]{Jin}. 
Finally, the front face is commutative because $\pi_\sZ$ is a triangle functor. 
%
%
%
\end{proof}

\bigskip
{\small\texttt{
  \noindent        
  \begin{tabular}{lll}
    \textit{\textrm{Contact:}}
      & nathan.broomhead@plymouth.ac.uk, & r.coelhosimoes@lancaster.ac.uk,  \\
      & d.pauksztello@lancaster.ac.uk,          & jonathan.woolf@liverpool.ac.uk         
  \end{tabular}}}
  
\end{document}